\documentclass{amsart}
\usepackage{amsfonts}

\usepackage{graphicx}
\usepackage{amscd}
\usepackage{amsmath}
\usepackage{amssymb}

\newtheorem{theorem}{Theorem}[section]
\theoremstyle{plain}

\newtheorem{claim}[theorem]{Claim}
\newtheorem{fact}[theorem]{Fact}

\newtheorem{definition}[theorem]{Definition}

\newtheorem{lemma}[theorem]{Lemma}

\newtheorem{question}{Question}

\newtheorem{proposition}[theorem]{Proposition}
\newtheorem{remark}[theorem]{Remark}

\numberwithin{equation}{section}

\def\eps{{\varepsilon}}

\def\dim{{\rm dim}}
\def\newEE{{Q}}

\begin{document}

\title[$K$-divisibility constants for some special couples]{On the $K$-divisibility constant for
some special finite-dimensional Banach couples}

\address{Yacin Ameur, Department of Chemistry and Biomedical Sciences, University of Kalmar, SE-391 82 Kalmar,
Sweden}
\email{Yacin.Ameur@hik.se}
\address{Michael Cwikel, Department of Mathematics, Technion - Israel Institute of
Technology, Haifa 32000, Israel}
\email{mcwikel@math.technion.ac.il}
\author[Yacin Ameur]{Yacin Ameur}
\author[Michael Cwikel]{Michael Cwikel}

\begin{abstract} We prove new estimates of the $K$-divisibility
constants for some special Banach couples. In particular, we prove
that the $K$-divisibility constant for a couple of the form
$(U\oplus V, U)$ where $U$ and $V$ are non-trivial Hilbert spaces
equals $2/\sqrt{3}$. We also prove estimates for the
$K$-divisibility constant of the two-dimensional version of the
couple $(L_2,L_\infty)$, proving in particular that this couple is
not exactly $K$-divisible. There are also several auxiliary
results, including some estimates for relative Calder\'on
constants for finite dimensional couples.
\end{abstract}

\maketitle \tableofcontents

\section{\label{intro}Introduction}
\noindent Let us begin by recalling the celebrated
Brudnyi-Krugljak $K$-divisibility theorem (cf.\ \cite{bk1},
\cite[p. 325, Paragraph C and Theorem 3.2.7 ]{bk}).

\begin{theorem}
\label{bkthm}Let $\vec{A}=(A_{0},A_{1})$ be a Banach couple, and let
$N$ be either a fixed natural number or $\infty$. There exists a
constant $C_N$, depending only on $\vec{A}$ and $N$, which has the
following property: Suppose that
$a$\ is an arbitrary element of $A_{0}+A_{1}$ whose Peetre $K$%
-functional satisfies the estimate
\begin{equation}
K(t,a;\vec{A})\le \sum_{n=1}^{N }\phi _{n}(t)\ \text{for all }t>0,
\label{plomj}
\end{equation}
where the functions $\phi _{n}$\ are each positive and concave on
$(0,\infty )$ and $\sum_{n=1}^{N }\phi _{n}(1)<\infty $. Then
there exists a sequence of elements $a_{n}\in A_{0}+A_{1}$\ such
that $a=\sum_{n=1}^{N }a_{n}$\ (where this series converges in
$A_{0}+A_{1}$ norm) and
\begin{equation}
K(t,a_{n};\vec{A})\le C_N\phi _{n}(t)\ \text{for all }t>0\text{
and each }1\le n< N+1.  \label{zronlk}
\end{equation}
\end{theorem}

The main interest of Theorem \ref{bkthm} resides in the special case
when $N=\infty$, but we will also need to consider other values of
$N$
below.
We refer to \cite{bk} and also to remarks in the
introductions of \cite {ckdiv} and \cite{cjm} for more details about
Theorem \ref{bkthm} and its applications. Its original proof appears
in \cite{bk}. Various alternative proofs using the so-called
``strong fundamental lemma'' can be found in \cite {ckdiv},
\cite{bsh}, cf.\ also \cite{bs1}.

We shall use the notation $\gamma_N (\vec{A})$ for the infimum of all numbers $%
C_N$ having the property stated in Theorem \ref{bkthm}. This
number may be called the \textit{$N$-term $K$-divisibility
constant for }$\vec{A}$\textit{\ }. When $N=\infty$, we follow the
notation and terminology of previous papers and simply write
$\gamma(\vec{A})$ instead of $\gamma_N(\vec{A})$ and speak of the
{\it $K$-divisibility constant } of $\vec{A}$. It is not hard to
check that these constants satisfy
$$1\le\gamma_i(\vec{A})\le \gamma_j(\vec{A})\le
\gamma({\vec{A}}),\quad 1\le i\le j.$$ (Strictly speaking, the first
inequality is only true if  $\vec{A}$ is {\it non-zero}. ``Non zero"
means that we exclude the trivial cases where $A_0=A_1$ and this
space consists solely of the zero element of some Hausdorff
topological vector space. In these cases $\gamma(\vec{A})=0$.)

All Banach spaces in this paper will be assumed to be over the
reals,
except
when it is explicitly stated otherwise. But it is
clear from the statement of Theorem \ref{bkthm} that if $A_{0}$
and $A_{1}$ happen to be complex Banach spaces, then the value of
$\gamma (\vec{A})$ will be the same, independently of whether we
consider the underlying scalar field to be $\Bbb{R}$ or $\Bbb{C}$.
For a related comment see Remark \ref{rc}.

Our main goal in this paper is to calculate the exact value of,
and obtain new estimates for $\gamma (\vec{A})$ for some
particular ``natural'' choices of the couple $\vec{A}$. Some of
the auxiliary results which we obtain en route to this goal may
perhaps also be useful in the future for other purposes, including
the determination of $\gamma (\vec{A})$ for other couples.

Theorem \ref{bkthm} is one of the most important and useful results in real
interpolation theory, and potentially also has interesting applications
beyond that theory. In its applications so far, the precise value of $\gamma
(\vec{A})$ does not seem to play a crucial role. However, as has turned out
to be the case with other important theorems in analysis, we believe that
searching for optimal constants, and thus optimal proofs, can also enhance
our general understanding of this very significant result.

It is known (cf.\ \cite{cjm}) that
\begin{equation}
1\le\gamma (\vec{A})\le 3+2\sqrt{2} \approx 5.\,8284 \ \ .
\label{cjmbe}
\end{equation}
for every non zero Banach couple $\vec{A}$.

Recently \cite{cw1} it has been shown that, in the case where $\vec{A}$ is a
non zero couple of Banach lattices (or complexified Banach lattices) of
measurable functions on the same underlying measure space, the estimate (\ref
{cjmbe}) can be sharpened to
\begin{equation*}
1\le \gamma (\vec{A})\le 4 .
\end{equation*}

A number of couples $\vec{A}$ are known to be \textit{exactly $K$-divisible}%
, i.e.\ to have the property that $\gamma(\vec{A})=1$. These include $%
(L^1,L^\infty)$ and the ``weighted" $L^p$ couples $(L^1_{w_0},L^1_{w_1})$
and $(L^\infty_{w_0},L^\infty_{w_1})$, for all choices of weight functions $%
w_0$ and $w_1$. The proof that $\gamma(\vec{A})=1$ for the first
of these couples can be found in \cite{daa}. It also follows from
an obvious generalization of the proof of Lemma 5.2 of \cite{cp}
p.\ 44. The proof for the latter two couples is contained in
Proposition 3.2.13 of \cite{bk} p.\ 335. Let us also mention
another collection of trivial examples of
exactly $K$-divisible couples. These are the non zero couples $\vec{A}%
=(A_0,A_1)$ for which $A_0=A_1$ isometrically. (For such a couple, every
element $a\in A_0+A_1$ satisfies $K(t,a;\vec{A})=\min\{1,t\}\|a\|_{A_0}$.
So, if $a$ satisfies (\ref{plomj}) and we choose $a_n=\frac{\phi _{n}(1)}{%
\sum_{m=1}^{\infty }\phi _{m}(1)}a$ for each $n \in \Bbb{N}$, then
it is obvious that we obtain (\ref{zronlk}) with $C_\infty=1$ when
$t=1$, and consequently also for all $t>0$.)

On the other hand it is also known that $\gamma(\vec{A})>1$ for certain
couples $\vec{A}$. The first example to be given of such a couple was $\vec{A%
}=(C,C^1)$, studied by Krugljak in \cite{kr}. Subsequently Podgova \cite{po}
showed that this same couple satisfies $\gamma (\vec{A}) \ge \frac{3+2\sqrt{2%
}}{1+2\sqrt{2}} \approx 1.\,5224 $. As announced in \cite{shv}, Pavel
Shvartsman has produced a different and much simpler example of a couple $%
\vec{S}=\left(S_0,S_1\right)$
whose
$2$-term $K$-divisibility constant
satisfies
$\gamma_2 (\vec{S}) = \frac{3+2\sqrt{2%
}}{1+2\sqrt{2}}$. He takes $S_{0}$ to be $\Bbb{R}^{2} $ equipped with the $%
\ell ^{\infty }$ norm and $S_{1}$ to be a one dimensional subspace of $\Bbb{R%
}^{2}$ whose unit ball is a line segment which makes an angle of
$\frac{\pi }{8}$ with one of the coordinate axes. Furthermore,
Shvartsman shows that this couple is ``extremal" among all couples
$\vec{A}=\left(A_0,A_1\right)$ satisfying $A_{j}\subset \Bbb{R}^{2}$
for $j=0,1$, in the sense that all such couples satisfy $\gamma_2
(\vec{A})\le \frac{3+2\sqrt{2}}{1+2\sqrt{2}}$. It will follow from
one of
our results in this paper that $\gamma(\vec{S})\le
2\sqrt{2/3}\approx 1.6330$, and thus that the
exact
value of
$\gamma(\vec{S})$
lies somewhere
in the interval $(1.52,1.64)$.

Apparently, neither $(C,C^1)$ nor Shvartsman's finite dimensional couple can
be realized as couples of Banach lattices on a measure space. But it turns
out that there also exist couples of lattices whose $K$-divisibility
constant is bigger than $1$. The first examples of such couples were found
in \cite{ckeich}. They are somewhat ``exotic" couples of spaces $\vec{A}%
=\left(A_0,A_1\right)$ both contained in $\Bbb{R}^3$. They each satisfy $%
\gamma(\vec{A})>1$ as a consequence of the fact that they do not possess
another property, \textit{almost exact monotonicity}, which is defined on
p.\ 30 of \cite{ckeich}.

In this paper we deal with what could be considered two of the
simplest, ``nicest" and most ``natural" couples among those which
are not already known to be exactly $K$-divisible, namely a couple
$\vec{H}=\left(H_0,H_1\right)$ of Hilbert spaces, and the lattice
couple $(L^2,L^\infty)$. In addition to its other good properties,
$\left(H_0,H_1\right)$ is known, as shown in \cite{ameur}, to be
an exact Calder\'on couple. $(L^2,L^\infty)$ is also a Calder\'on
couple \cite{lorshi} and the optimal decomposition for obtaining
its $K$-functional exactly is quite simple to describe. But it
turns out, perhaps rather surprisingly,
that neither of these couples are exactly $K$%
-divisible in general, and one can even find two-dimensional
versions of each of these couples for which
exact $K$-divisibility does not hold.

The paper is organized as follows: In Section \ref{prelim} we recall some
definitions and collect some general preliminary results which will be
needed in other sections.
In Section \ref{hc} we find the exact value of $%
\gamma (\vec{Y})$ where $\vec{Y}$ is the simplest non trivial
version of a couple of Hilbert spaces. Our result is that
$\gamma(\vec{Y})=2/\sqrt{3}$. After considering various
generalizations of this result, we consider all other couples of
(real) Hilbert spaces which are contained in $\Bbb{R}^{2}$, and we
prove a (rather more crude) upper estimate for their
$K$-divisibility constants, namely $\gamma(\vec{G})<\sqrt{2}$.

Finally, in Section \ref{latticecouple} we consider the couple
$(L^{2},L^{\infty })$ and, in particular, the case where the
underlying measure space consists of two atoms of equal measure,
i.e.\ the two dimensional couple $\vec{X}=(\ell _{2}^{2},\ell
_{2}^{\infty }) $. It turns out to be quite easy to show that
$\gamma (\vec{X})>1$. But the determination of the exact value of
$\gamma (\vec{X})$ is a much longer and as yet unfinished story. We
obtain some (rather complicated) equations
which in principle could be solved to obtain the exact value of $\gamma (%
\vec{X})$. Numerical experiments suggest that maybe $\gamma
(\vec{X})$ is approximately equal to $1.03$. The sharpest
estimates which we have are
$$1<\gamma(\vec{X})<\frac{4+3\sqrt{2}}{4+2\sqrt{2}}\approx 1.\,2071\ \ .$$

In the Appendix we prove that the couple $(L^2,L^\infty)$ is an
exact Calder\'{o}n couple in the two-dimensional case, but not in
the eight-dimensional case. This example proves that there
are
in general no tight connections between the properties of
being an exactly $K$-divisible couple and of being an exact Calder\'on couple.

\section{\label{prelim}Some definitions and general preliminary results}

For the basic notions of the real method of interpolation, we
refer, e.g.\
to \cite{bsh}, \cite{bl} or \cite{bk}. For any given Banach couple $\vec{A}%
=\left( A_{0},A_{1}\right) $, we let $A_{j}^{\sim }$ denote the \textit{%
Gagliardo completion }of $A_{j}$, $j=0,1$, i.e.\ the Banach space
of elements $a$ of $A_{0}+A_{1}$ which are limits in $A_{0}+A_{1}$
norm of bounded sequences in $A_{j}$ or, equivalently, for which
the norm $\Vert a\Vert _{A_{j}^{\sim
}}=\sup_{t>0}K(t,a;\vec{A})/t^{j}$ is finite. Obviously
$A_{0}^{\sim }+A_{1}^{\sim }=A_{0}+A_{1}$. We also recall that the couple $%
\vec{A}=\left( A_{0},A_{1}\right) $ and the corresponding couple
of its Gagliardo completions $\vec{A^{\sim }}=\left( A_{0}^{\sim
},A_{1}^{\sim
}\right) $ have identical $K$-functionals, i.e.\ $K(t,a;\vec{A})=K(t,a;\vec{%
A^{\sim }})$ for all $a\in A_{0}+A_{1}$ and all $t>0$.
Consequently we also have $\gamma (\vec{A})=\gamma (\vec{A^{\sim
}})$.

There is a close connection between $K$-divisibility and couples
of weighted $L^{1}$ spaces which we wish to exploit. Our point of
departure is the following lemma.

\begin{lemma}
\label{sgk}\smallskip Let $\vec{A}=(A_{0},A_{1})$ be an arbitrary
Banach couple and let $a$ be an arbitrary element of
$A_{0}+A_{1}$. Then there
exist a measure space $(\Omega ,\mathcal{S},\mu )$ and measurable functions $%
w_{j}:\Omega \rightarrow (0,\infty ]$ for $j=0,1$ and a measurable function $%
f_{a}:\Omega \rightarrow [0,\infty )$ such that $K(t,a;\vec{A})=K(t,f_{a};%
\vec{P})$ for all $t>0$, where $\vec{P}$ is the couple of weighted
$L^{1}$ spaces $\vec{P}=(L_{w_{0}}^{1}(\mu ),L_{w_{1}}^{1}(\mu
))$.
\end{lemma}

The straightforward proof of this result, which uses \cite[Lemma
5.4.3, p.\ 117]{bl}, can be found in \cite[pp.\ 46--47]{ckdiv}. It
should not be overlooked that the weight functions $w_{0}$ and
$w_{1}$ in Lemma \ref{sgk} have the slightly exotic property that
they are permitted to assume the value $+\infty $. Since every
function in $L_{w_{0}}^{1}(\mu )+L_{w_{1}}^{1}(\mu )$ vanishes
a.e.\ on the set where $w_{0}=w_{1}=\infty $, we always can and
will assume that this set is empty. We also mention that the proof
in \cite{ckdiv} shows that $(\Omega ,\mathcal{S},\mu )$ and
$w_{0}$ and $w_{1}$ can be chosen rather simply and quite
explicitly, and we can also, for example, arrange things so that
$f_a$ is a constant function.

It turns out that for each $\vec{A}$ and each $a\in A_{0}+A_{1} $
and each corresponding $\vec{P}$ and $f_{a}$ with the properties
just specified, there exists a bounded linear operator
$T:\vec{P}\rightarrow \vec{A^{\sim }}$ such that $a=Tf_{a}$. Let
$\mathcal{T}_{a}$ denote the set of all such operators $T$ for some
given choice of $a$ and $f_{a}$. Then it turns out that
\begin{equation}
\gamma (\vec{A})=\sup_{a\in A_{0}+A_{1}} c_a \quad \text{where}\quad
c_a=c_a(\vec{A}):= \inf_{T\in \mathcal{T}%
_{a}}\left\| T\right\| _{\vec{P}\rightarrow \vec{A^{\sim }}}.
\label{rbb}
\end{equation}

This formula, whose proof will be briefly recalled below, turns
out to be particularly suitable for our calculations of
$K$-divisibility constants in this paper.

It is sometimes convenient to re-express (\ref{rbb}) slightly
differently. For $\vec{A}$, $a$, $\vec{P}$ and $f_{a}$ as above,
let $\Lambda _{a}$ be
the set of linear operators $T:\vec{P}\rightarrow \vec{A^{\sim }}$ with $%
\left\| T\right\| _{\vec{P}\rightarrow \vec{A^{\sim }}}\le 1$ such that $%
Tf_{a}=\lambda a$ for some positive number $\lambda =\lambda
_{T}$. Then obviously (\ref{rbb}) is the same as
\begin{equation}
\gamma (\vec{A})=\sup_{a\in A_{0}+A_{1}}\left( \inf_{T\in \Lambda _{a}}\frac{%
1}{\lambda _{T}}\right) .  \label{ricv}
\end{equation}

\begin{remark}\smallskip \label{4lat}
Clearly $\mathcal{T} _{ta}=\mathcal{T} _{a}$ and so $c_{ta}=c_a$ for
all scalars $t\ne 0$. Furthermore, if, as is the case for most
couples considered in the paper, $A_{0}$ and $A_{1}$ are both Banach
lattices of measurable functions on the same underlying measure
space, then it is easy to see that, in the formula (\ref{rbb}), the
supremum can be replaced by the supremum over all non negative
functions $a$ in $A_{0}+A_{1}$.

Indeed, we have for every $a\in A_0+A_1$ that $c_a=c_{|a|}$.
\end{remark}

\smallskip

At first sight it seems that there could be some ambiguity in
(\ref{rbb}), because the set $\mathcal{T}_{a}$ depends on our
particular choices of the
measure space $(\Omega ,\mathcal{S},\mu )$ and the associated functions $%
f_{a}$, $w_{0}$ and $w_{1}$. The key to showing that in fact there
is no such ambiguity is the theorem of Sedaev-Semenov
\cite{sedsem} (see \cite {cwikoz} for an alternative proof) or,
more precisely, the generalization of that theorem \cite[Theorem
3, pp.\ 47--49]{ckdiv} to the case of weight functions which are
permitted to take the value $+\infty $. It follows immediately
from that theorem, that if $(\Xi ,\mathcal{Y},\sigma )$ is a
second measure space and $v_{0}$ and $v_{1}$ are weight functions
and $g_{a}$
is a non negative measurable function such that $K(t,g_{a};L_{v_{0}}^{1}(%
\sigma ),L_{v_{1}}^{1}(\sigma ))=K(t,f_{a};L_{w_{0}}^{1}(\mu
),L_{w_{1}}^{1}(\mu ))$ for all $t>0$ then, for each $\epsilon
>0$, there exist two linear operators $U:(L_{v_{0}}^{1}(\sigma
),L_{v_{1}}^{1}(\sigma
))\rightarrow (L_{w_{0}}^{1}(\mu ),L_{w_{1}}^{1}(\mu ))$ and $%
V:(L_{w_{0}}^{1}(\mu ),L_{w_{1}}^{1}(\mu ))\rightarrow
(L_{v_{0}}^{1}(\sigma
),L_{v_{1}}^{1}(\sigma ))$ which satisfy $Ug_{a}=f_{a}$ , $Vf_{a}=g_{a}$, $%
\left\| U\right\| _{(L_{v_{0}}^{1}(\sigma ),L_{v_{1}}^{1}(\sigma
))\rightarrow (L_{w_{0}}^{1}(\mu ),L_{w_{1}}^{1}(\mu ))}\le
1+\epsilon $ and $\left\| V\right\| _{(L_{w_{0}}^{1}(\mu
),L_{w_{1}}^{1}(\mu ))\rightarrow (L_{v_{0}}^{1}(\sigma
),L_{v_{1}}^{1}(\sigma ))}\le 1+\epsilon $ .

\smallskip By composing the operators $U$ and $V$ with other suitable
operators, we readily see that the quantity $\inf_{T\in \mathcal{T}%
_{a}}\left\| T\right\| _{\vec{P}\rightarrow \vec{A^{\sim }}}$ is
independent
of the choices of the measure space, weight functions and the function $%
f_{a} $.

\smallskip For the convenience of the reader who may not be familiar with
these details, we mention that the fact that $\mathcal{T}_{a}$ is
non empty and the formula (\ref{rbb}) are both obtained by
considering the following theorem which, as we shall explain, is
intimately related, in fact equivalent, to Theorem \ref{bkthm}
(Cf.\ \cite[Proposition 1.40]{cn}).

\def \cconst {M}

\begin{theorem}
\smallskip \label{ec}Let $\vec{A}=(A_{0},A_{1})$ be an arbitrary Banach
couple. Then there exist constants $\cconst_{1}$, $\cconst_{2}$ and $\cconst_{3}$,
depending only on $\vec{A}$, with, respectively, the following
properties:

(i) For each $a\in A_{0}+A_{1}$, there exists a sequence $\left\{
a_{\nu }\right\} _{\nu \in \Bbb{Z}}$ of elements in $A_{0}^{\sim
}\cup A_{1}^{\sim }
$ which satisfies $a=\sum_{\nu \in \Bbb{Z}}a_{\nu }$ (convergence in $%
A_{0}+A_{1}$ norm) and also
\begin{equation}
\sum_{\nu \in \Bbb{Z}}\min \left\{ \left\| a_{\nu }\right\|
_{A_{1}^{\sim
}},t\left\| a_{\nu }\right\| _{A_{1}^{\sim }}\right\} \le \cconst_{1}K(t,a;\vec{A})%
\text{ for all }t>0.  \label{dvap}
\end{equation}

(ii) Let $w_{0}$ and $w_{1}$ be arbitrary weight functions on an
arbitary measure space $(\Omega ,\mathcal{S},\mu )$. Let $\vec{P}$
be the couple of weighted $L^{1}$ spaces
$\vec{P}=(L_{w_{0}}^{1}(\mu ),L_{w_{1}}^{1}(\mu ))$. Suppose that
the elements $a\in A_{0}+A_{1}$ and $f\in
L_{w_{0}}^{1}+L_{w_{1}}^{1}$ satisfy
\begin{equation}
K(t,a;\vec{A})\le K(t,f;\vec{P})\text{ for all }t>0.  \label{rbv}
\end{equation}
Then there exists a bounded linear operator $T:\vec{P}\rightarrow \vec{%
A^{\sim }}$ such that $\left\| T\right\| _{\vec{P}\rightarrow \vec{A^{\sim }}%
}\le \cconst_{2}$ and $Tf=a$.

(iii) Suppose that $(\Omega ,\mathcal{S},\mu )$, $w_{0}$, $w_{1}$,
$f$ and $a $ are exactly as in part (ii), except that instead of
(\ref{rbv}) they satisfy
\begin{equation*}
K(t,a;\vec{A})\mathbf{=}K(t,f;\vec{P})\text{ for all }t>0.
\end{equation*}
Then there exists a bounded linear operator $T:\vec{P}\rightarrow \vec{%
A^{\sim }}$ such that $\left\| T\right\| _{\vec{P}\rightarrow \vec{A^{\sim }}%
}\le \cconst_{3}$ and $Tf=a$.

In fact the infima of all constants $\cconst_1$, $\cconst_2$ and
$\cconst_3$ satisfying (i), (ii) and (iii) respectively, coincide,
and they all equal $\gamma(\vec{A})$, the infimum of the constants
$C_\infty$ for which Theorem \ref{bkthm} holds.
\end{theorem}

For a proof of part (ii) of this theorem, which uses Theorem
\ref{bkthm} and gives the value $\cconst_{2}=C_\infty+\epsilon $ for any
choice of $\epsilon >0$, see \cite[Theorem 4.4.12, pp.\
586--588]{bk}. We mention in passing that part (ii) has an
important and immediate consequence. It provides a simple
description of all relative interpolation spaces for operators
mapping from any weighted $L^{1}$ couple into any Banach couple
$\vec{A}$ which satisfies $A_{j}^{\sim }=A_{j}$ for $j=0,1$.

\smallskip Part (i), also known as the ``strong fundamental lemma'', is
proved in \cite[Theorem 4, pp.\ 59--54]{ckdiv} for
$\cconst_{1}\approx 8$ and, with a better constant
$\cconst_{1}\approx 3+2\sqrt{2}$, in \cite[pp.\ 73--77]{cjm}. Cf.\
also \cite{cw1} for more explicit versions of some of the steps of
the proof in \cite{cjm}. (Note that in (\ref{dvap}) we adopt the
conventions that $\left\| a_{\nu }\right\| _{A_{j}^{\sim }}=\infty $
if $a\notin A_{j}^{\sim }$ and that $\min \{\alpha ,\infty \}=\min
\{\infty ,\alpha \}=\alpha $ for every $\alpha \in \Bbb{R}$.)

Part (ii) can be deduced from part (i), and with
$\cconst_{2}=\cconst_{1}+\epsilon $ for any choice of $\epsilon >0$. This can be
done, using (an obvious modification of) an argument which appears
in \cite[pp.\ 54--55]{ckdiv} cf.\ also \cite[Theorem 4.8, p.\
38]{cp}. Moreover, this result, and also part (iii), are also both
valid in the case where either or both of the weight functions
$w_{0}$ and $w_{1}$ are permitted to take the value $+\infty $ on
some subsets of $\Omega $. The proof in \cite{ckdiv} makes use of
the generalized version \cite[Theorem 3, p.\ 47]{ckdiv} of the
Sedaev-Semenov theorem already mentioned above. (The Sedaev-Semenov
theorem is also the main, perhaps only, ingredient of the ``obvious
modification'' mentioned above.)

The connection between parts (ii) and (iii) is a simple matter.
Obviously (ii) implies (iii) with $\cconst_{3}\mathbf{=}\cconst_{2}$. On the
other hand we can also easily obtain that (iii) implies (ii) with
$\cconst_{2}=\cconst_{3}+\epsilon $ for any choice of $\epsilon >0$. This is
done by first using Lemma \ref{sgk} to obtain $f_{a}$ and then
using the generalized version of the Sedaev-Semenov
theorem to find a linear map $U$ between appropriate couples of weighted $%
L^{1}$ spaces, which satisfies $Uf=f_{a}$ and has norm arbitrarily close to $%
1$.

Theorem \ref{bkthm}, with $C_\infty=\cconst_{2}$ can be deduced from
part (ii) of Theorem \ref{ec}, again using arguments from
\cite[pp.\ 54--55]{ckdiv} and using the more general version where
the weight functions are permitted to take infinite values.

Conversely, as mentioned in \cite[p.\ 71]{cjm} and shown more
explicitly in \cite[Proposition 1.40]{cn}, it is also possible to
deduce part (i) (and consequently also part (ii)) of Theorem
\ref{ec} from Theorem \ref{bkthm}, with $\cconst_{1}=C_\infty+\epsilon $
for any choice of $\epsilon >0$.

\smallskip

It should be noted that part (iii) of the above theorem, together
with the connections described above between the constants
$\cconst_{1}$, $\cconst_{2}$ and $\cconst_{3}$ for which parts
(i), (ii) and (iii) of the theorem hold, give us the formula
(\ref{rbb}).

\smallskip

For most couples $\vec{A}=(A_{0},A_{1})$ which we study in this
paper, $A_{0}$ and $A_{1}$ are both finite dimensional. For such
couples it is clear that $A_{j}^{\sim }=A_{j}$ isometrically for
$j=0,1$. It is also helpful to know, as the following lemma shows,
that, for such couples, the infimum $\inf_{T\in
\mathcal{T}_{a}}\left\| T\right\| _{\vec{P}\rightarrow \vec{A}}$
appearing in (\ref{rbb}) is actually attained for each fixed
element $a$. This of course implies that the infimum $\inf_{T\in
\Lambda _{a}}1/\lambda _{T}$ in (\ref{ricv}) is also attained for
each $a$. We will refer to any operator $T$ for which this latter
infimum is attained as an \textit{\textbf{optimal element of}}
$\Lambda _{a}$. Obviously such an operator satisfies $\left\|
T\right\| _{\vec{P}\rightarrow \vec{A}}=1$.

\begin{lemma}
\label{findim}Let $\vec{F}=(F_{0},F_{1})$ and
$\vec{A}=(A_{0},A_{1})$ be Banach couples and suppose that
$A_{0}+A_{1}$ is a finite dimensional space. Let $a$ and $f$ be
arbitrary fixed elements of $A_{0}+A_{1}$ and $F_{0}+F_{1} $
respectively. Suppose that the class $\mathcal{T}_{a}$ of all
bounded linear operators $T:\vec{F}\rightarrow \vec{A}$ which
satisfy $Tf=a$ is non
empty. Then there exists an operator $S\in \mathcal{T}_{a}$ such that $%
\left\| S\right\| _{\vec{F}\rightarrow \vec{A}}=\inf_{T\in \mathcal{T}%
_{a}}\left\| T\right\| _{\vec{F}\rightarrow \vec{A}}$.
\end{lemma}

\textit{Proof.} Let $N$ be the dimension of $A_{0}+A_{1}$ and let $%
\{e_{k}\}_{k=1}^{N}$ be a basis of $A_{0}+A_{1}$. Then every
bounded operator $T:F_{0}+F_{1}\rightarrow A_{0}+A_{1}$ defines
and can be defined by a collection $\lambda _{1},\lambda
_{2},...,\lambda _{N}$ of $N$ linear
bounded linear functionals on $F_{0}+F_{1}$, via the formula $%
Tg=\sum_{k=1}^{N}\lambda _{k}(g)e_{k}$ for each $g\in
F_{0}+F_{1}$. Consider
a sequence of elements $\left\{ T_{n}\right\} _{n\in \Bbb{N}}$ in $\mathcal{T%
}_{a}$ such that $\left\| T_{n}\right\| _{\vec{F}\rightarrow
\vec{X}}\le
c_{a}+1/n$, where $c_{a}=\inf_{T\in \mathcal{T}_{a}}\left\| T\right\| _{\vec{%
F}\rightarrow \vec{A}}$. Let $\lambda _{n,k}$ denote the bounded
linear functional on $F_{0}+F_{1}$ defined for each $n\in \Bbb{N}$
and each $k\in \{1,2,...,N\}$, such that
$T_{n}g=\sum_{k=1}^{N}\lambda _{n,k}(g)e_{k}$ for each $g\in
F_{0}+F_{1}$. Now let us define the operator $S$ by
\begin{equation*}
Sg=\sum_{k=1}^{N}\lambda _{*,k}(g)e_{k}\text{ for each }g\in
F_{0}+F_{1},
\end{equation*}
where the $N$ linear functionals $\lambda _{*,1},\lambda
_{*,2},...,\lambda _{*,N}$ and $\lambda _{1}$ are given by
$\lambda _{*,k}(g)=B\left( \left\{ \lambda _{n,k}(g)\right\}
_{n\in \Bbb{N}}\right) $ for each $g\in F_{0}+F_{1} $, where $B\in
\left( \ell ^{\infty }\right) ^{*}$ is a Banach limit, (i.e.\ an
element of $\left( \ell ^{\infty }\right) ^{*}$ which satisfies
$\left| B\left( \left\{ u_{n}\right\} _{n\in \Bbb{N}}\right)
\right| \le \limsup_{n\rightarrow \infty }|u_{n}|$ for all
$\left\{ u_{n}\right\} _{n\in \Bbb{N}}\in \ell ^{\infty }$ and
also $B\left( \left\{ u_{n}\right\} _{n\in \Bbb{N}}\right)
=\lim_{n\rightarrow \infty }u_{n}$ for every convergent sequence
$\left\{ u_{n}\right\} _{n\in \Bbb{N}}$ ). It easy to see that
each sequence $\left\{ \lambda _{n,k}(g)\right\} _{n\in \Bbb{N}}$
is indeed in $\ell ^{\infty }$ and it is straightforward, if a
little tedious, to verify that the operator $S$ has all the
required properties. We leave these matters to the reader.  \qed

\smallskip

\begin{remark}
\smallskip \label{rc}For all the couples $\vec{A}=(A_{0},A_{1})$
considered in this paper, $A_{0}$ and $A_{1}$ are both Banach
lattices of real valued measurable functions with the same
underlying measure space. As usual, we can define the
complexification of such a lattice $A_{j}$ to be the space, which
we may denote by $A_{j}^{\Bbb{C}}$, consisting of all
complex valued measurable functions $g$ such that $\left| g\right| \in A_{j}$%
, with the obvious norm. It is easy to see that the complexified
lattice
couple $\vec{A^{\Bbb{C}}}=(A_{0}^{\Bbb{C}},A_{1}^{\Bbb{C}})$ satisfies $%
\gamma (\vec{A^{\Bbb{C}}})\mathbf{=}\gamma (\vec{A})$. (Use the
fact that for any function $a\in A_{0}^{\Bbb{C}}+A_{1}^{\Bbb{C}}$
we have $K(t,a;\vec{A^{\Bbb{C}}})=K(t,|a|;\vec{A})$.)
\end{remark}

\section{\label{hc}On the $K$-divisibility constant for Hilbert
couples}
\subsection{\label{hc1}The $K$-divisibility constant for the couple $\vec{Y}%
=\left( Y_{0},Y_{1}\right) =\left( \ell _{2}^{2},\ell _{1}^2
\right) $.}

\smallskip The purpose of this subsection is to prove the following
theorem.

\begin{theorem}
\label{cfy}Let $\vec{Y}=\left( Y_{0},Y_{1}\right) $ be the Banach couple of
subspaces of $\Bbb{R}^{2}$ obtained by taking the unit ball of $Y_{0}$ to be
the disk $\{(x,y)\in \Bbb{R}^{2}:x^{2}+y^{2}\le 1\}$ and the unit ball of $%
Y_{1}$ to be the line segment $\left\{ (x,0)\in \Bbb{R}^{2}:-1\le x\le
1\right\} $. Then the $K$-divisibility constant of the couple $\vec{Y}$ is
given by
\begin{equation}
\gamma (\vec{Y})=\frac{2}{\sqrt{3}}.  \label{cy}
\end{equation}
\end{theorem}

\textit{Proof.} Consider the point $\alpha =(\cos a,\sin a)\in Y_{0}+Y_{1}$
where $a\in [0,2\pi )$. Let $E_{a}$ be the set consisting of every number
which is the norm $\left\| T\right\| _{\vec{P}\rightarrow \vec{Y}}$ of some
bounded linear operator $T$ from some couple $\vec{P}$ of weighted $L^{1}$
spaces into $\vec{Y}$, which satisfies $Tf=\alpha $ for some element $f\in
P_{0}+P_{1}$ for which
\begin{equation}
K(t,f;\vec{P})=K(t,\alpha ;\vec{Y})\text{ for all }t>0.  \label{kfalf}
\end{equation}
Note that the weight functions $w_{0}$ and $w_{1}$ used in the
definition of $P_{0}$ and $P_{1}$ are permitted to assume the
value $+\infty $ on some sets of positive measure. We shall
explicitly need this option here.

Let $c_{a}=\inf E_{a}$. It follows from Remark \ref{4lat} that
$\gamma(\vec{Y})=\sup_{a\in[0,\pi/2]}c_a$. We claim that in fact
\begin{equation}
\gamma (\vec{Y})=\sup_{a\in (0,\pi /2)}c_{a}.  \label{opin}
\end{equation}

To show (\ref{opin}) we first observe that, since $K(t,Tf;\vec{Y})\le
\left\| T\right\| _{\vec{P}\rightarrow \vec{Y}}K(t,f;\vec{P})$ for all $t>0$
and for every bounded operator $T:\vec{P}\rightarrow \vec{Y}$, we must have $%
c_{a}\ge 1$ for every $a\in [0,\pi /2]$. It turns out to be rather easy to
show that $c_{a}\le 1$ in the two special cases, $a=0$ and $a=\pi /2$, and
this will of course imply (\ref{opin}).

In the case where $a=0$, i.e., $\alpha =(1,0)$, we use a very simple couple $%
\vec{P}$ where the underlying measure space consists of a single point $b$
which has measure $1$ and $\left\| h\right\| _{P_{0}}=\left\| h\right\|
_{P_{1}}=\left| h(b)\right| $ for every $h\in P_{0}+P_{1}$. (I.e., $%
w_{0}(b)=w_{1}(b)=1$.) We also use the ``function'' $f\in P_{0}+P_{1}$
defined by $f(b)=1$ which clearly satisfies
\begin{equation*}
K(t,f;\vec{P})=\min \{1,t\}=K(t,(1,0);\vec{Y})\text{ for all }t>0.
\end{equation*}
Then we use the operator $T$ defined by $T(h)=\left( h(b),0\right) $ for all
$h\in P_{0}+P_{1}$ to show that $c_{0}\le 1$.

In the case where $a=\pi /2$, i.e., $\alpha =(0,1)$, it is convenient, once
again, to use an underlying measure $(\Omega ,\Sigma ,\mu )$ space
containing (at least) one point $b$ which is an atom of measure $1$. But
this time the weight functions $w_{0}$ and $w_{1}$ for which $%
P_{j}=L_{w_{j}}^{1}(\mu )$ should be chosen to satisfy $w_{0}(b)=1$ and $%
w_{1}(b)=+\infty $. This means that every function $h$ in $P_{1}$ satisfies $%
h(b)=0$ and so the linear map $T$ defined by $Th=\left( 0,h(b)\right) $ maps
$P_{1}$ into $Y_{1}$ with norm $0$ and $P_{0}$ into $Y_{0}$ with norm $1$.
Furthermore the function $f=\chi _{\{b\}}$ satisfies $Tf=(0,1)$ and $K(t,f;%
\vec{P})=1=K(t,(0,1);\vec{Y})$ for all $t>0$. This shows that $c_{\pi /2}\le
1$ and so completes the proof of (\ref{opin}).

In the light of the preceding calculations it remains to calculate or
estimate $c_{a}$ for values of $a\in (0,\pi /2)$. So let us indeed fix $a\in
(0,\pi /2)$ and set $\alpha =(\cos a,\sin a)=(\alpha _{1},\alpha _{2})$. It
is easy to see that the error functional $E(t,\alpha ;\vec{Y})=\inf \left\{
\left\| \alpha -\beta \right\| _{Y_{0}}:\beta \in Y_{1},\left\| \beta
\right\| _{Y_{1}}\le t\right\} $ is given by the formula
\begin{equation*}
E(t,\alpha ;\vec{Y})=\left\{
\begin{array}{lll}
\sqrt{(t-\alpha _{1})^{2}+\alpha _{2}^{2}} & , & t\in [0,\alpha _{1}] \\
\alpha _{2} & , & t>\alpha _{1}.
\end{array}
\right. .
\end{equation*}

\smallskip Now we will describe a particular couple of weighted $L^{1}$
spaces $\vec{P}=(P_{0},P_{1})$ on the (non empty) interval $[0,\alpha _{1}]$%
, for which the function $f=\chi _{[0,\alpha _{1}]}$ satisfies (\ref{kfalf}%
). Once again we use the fact that (\ref{kfalf}) is equivalent to
\begin{equation}
E(t,f;\vec{P})=E(t,\alpha ,\vec{Y})\text{ for all }t>0.  \label{efalf}
\end{equation}
To make (\ref{efalf}) hold, we choose a measure $\mu $ on
$[0,\alpha _{1}]$ which coincides with Lebesgue measure on
$[0,\alpha _{1})$ and such that the singleton set $\left\{ \alpha
_{1}\right\} $ has measure $\mu \left( \left\{ \alpha _{1}\right\}
\right) =1$. Then we take $P_{j}=L_{w_{j}}([0,\alpha _{1}],\mu )$
for $j=0,1$, where the weight functions $w_{0}$ and $w_{1}$ are
defined by
\begin{equation*}
w_{0}(t)=\left\{
\begin{array}{lll}
-\frac{d}{dt}E(t,\alpha ;\vec{Y})=\frac{\alpha _{1}-t}{\sqrt{(\alpha
_{1}-t)^{2}+\alpha _{2}^{2}}} & , & t\in [0,\alpha _{1}) \\
\alpha _{2} & , & t=\alpha _{1}.
\end{array}
\right.
\end{equation*}
and
\begin{equation*}
w_{1}(t)=\left\{
\begin{array}{lll}
1 & , & t\in [0,\alpha _{1}) \\
+\infty & , & t=\alpha _{1}.
\end{array}
\right. .
\end{equation*}

Since $w_{0}$ is decreasing on $[0,\alpha _{1})$ it is easy to
obtain that
$$E(t,f,\vec{P})=\int_{[\min (t,\alpha _{1}),\alpha _{1}]}w_{0}d\mu
=\int_{[t,\infty )\cap [0,\alpha _{1})}w_{0}d\mu +\alpha _{2}$$ for each $%
t>0 $ which immediately also gives us (\ref{efalf}) and (\ref{kfalf}).

Let us now define $E_{a}^{*}$ to be the subset of $E_{a}$
consisting of the numbers $\left\| T\right\| _{\vec{P}\rightarrow
\vec{Y}}$ obtained in the special case where $\vec{P}$ is the
particular couple $$\vec{P}=( L_{w_{0}}^{1}([0,\alpha _{1}],\mu
),L_{w_{1}}^{1}([0,\alpha _{1}],\mu
))$$ which we have just defined, and the function $f$ for which $%
Tf=\alpha $ is given by $f=\chi _{[0,\alpha _{1}]}$. In view of
the generalized version of the Sedaev-Semenov theorem in \cite
{ckdiv}, it is clear that $c_{a}$ is also the infimum of the set
$E_{a}^{*}$.

Any bounded linear operator $T:\vec{P}\rightarrow \vec{Y}$ for this
particular choice of $\vec{P}$ must be given by the formula
\begin{equation}
Th=\left( \int_{[0,\alpha _{1})}g_{1}(\xi )h(\xi )d\xi +\beta _{1}h(\alpha
_{1}),\int_{[0,\alpha _{1})}g_{2}(\xi )h(\xi )d\xi +\beta _{2}h(\alpha
_{1})\right)  \label{deft}
\end{equation}
for all $h\in P_{0}+P_{1}$. Here $g_{1}$ and $g_{2}$ are suitable bounded
measurable functions on $[0,\alpha _{1})$ and $\beta _{1}$ and $\beta _{2}$
are real numbers. For all $h\in P_{1}$ we have $h(\alpha _{1})=0$. But all
such functions $h$ must also satisfy $\int_{[0,\alpha _{1})}g_{2}(\xi )h(\xi
)d\xi +\beta _{2}h(\alpha _{1})=0$. Consequently $g_{2}=0$ a.e.\ on $%
[0,\alpha _{1})$. Thus the norm $\left\| T\right\| _{P_{1}\rightarrow Y_{1}}$
equals $\left\| g_{1}\right\| _{L^{\infty }[0,\alpha _{1})}$. The norm $%
\left\| T\right\| _{P_{0}\rightarrow Y_{0}}$ is the supremum of
\begin{eqnarray*}
&&\theta _{1}\int_{[0.,\alpha _{1}]}\left( g_{1}\chi _{[0,\alpha
_{1})}+\beta _{1}\chi _{\{\alpha _{1}\}}\right) hd\mu +\theta _{2}\beta
_{2}h(\alpha _{1}) \\
&=&\int_{[0.,\alpha _{1}]}\left( \theta _{1}g_{1}\chi _{[0,\alpha
_{1})}+(\theta _{1}\beta _{1}+\theta _{2}\beta _{2})\chi _{\{\alpha
_{1}\}}\right) hd\mu
\end{eqnarray*}
as $h$ ranges over the unit ball of $P_{0}$ and $(\theta _{1},\theta _{2})$
ranges over the unit circle. Let us first calculate the supremum, for a
fixed choice of $(\theta _{1},\theta _{2})$, as $h$ ranges over the unit
ball of $P_{0}$. The standard duality between $L^{1}$ and $L^{\infty }$
gives us that this supremum equals
\begin{eqnarray}
&&\left\| \frac{\theta _{1}g_{1}\chi _{[0,\alpha _{1})}+(\theta _{1}\beta
_{1}+\theta _{2}\beta _{2})\chi _{\{\alpha _{1}\}}}{w_{0}}\right\|
_{L^{\infty }([0,\alpha _{1}],\mu )}  \notag \\
&=&\max \left\{ \theta _{1}\mathop{\mathrm{ess\ sup}}_{\xi \in [0,\alpha
_{1})}\left| \frac{g_{1}(\xi )}{w_{0}(\xi )}\right| ,\frac{\left| \theta
_{1}\beta _{1}+\theta _{2}\beta _{2}\right| }{\alpha _{2}}\right\} .
\label{gling}
\end{eqnarray}
We now claim that
\begin{equation}
\left\| T\right\| _{P_{0}\rightarrow Y_{0}}=\max \left\{ \mathop{\mathrm{ess%
\ sup}}_{\xi \in [0,\alpha _{1})}\left| \frac{g_{1}(\xi )}{w_{0}(\xi )}%
\right| ,\frac{\sqrt{\beta _{1}^{2}+\beta _{2}^{2}}}{\alpha _{2}}\right\} .
\label{ching}
\end{equation}
This is because the expression in (\ref{gling}) equals the
expression on the right side of (\ref{ching}) for a suitable
choice of $(\theta _{1},\theta _{2})$ on the unit circle (either
$(1,0)$ or $\left( \frac {\beta _{1}}{\sqrt{\beta
_{1}^{2}+\beta _{2}^{2}}},\frac {\beta _{2}}{\sqrt{\beta _{1}^{2}+\beta _{2}^{2}}}%
\right) $). Furthermore it is dominated by the expression on the right side
of (\ref{ching}) for all other points $(\theta _{1},\theta _{2})$ on the
unit circle.

Since $w_{0}(\xi )<1$ for all $\xi \in [0,\alpha _{1})$, we have that
\begin{eqnarray*}
\left\| T\right\| _{P_{1}\rightarrow Y_{1}} &=&\left\| g_{1}\right\|
_{L^{\infty }[0,\alpha _{1})} \\
&=&\mathop{\mathrm{ess\ sup}}_{[0,\alpha _{1})}\left| g_{1}(\xi )\right| \le %
\mathop{\mathrm{ess\ sup}}_{\xi \in [0,\alpha _{1})}\left| \frac{g_{1}(\xi )%
}{w_{0}(\xi )}\right| .
\end{eqnarray*}
This means that the norm $\left\| T\right\| _{\vec{P}\rightarrow \vec{Y}}$
is also given by the expression on the right side of (\ref{ching}).

Of course here we are only concerned with those operators $T$ for which $%
T\chi _{[0,\alpha _{1}]}=(\alpha _{1},\alpha _{2})$, i.e.\
\begin{equation}
\int_{\lbrack 0,\alpha _{1})}g_{1}(\xi )d\xi +\beta _{1}=\alpha _{1}\text{
and }\beta _{2}=\alpha _{2}.  \label{blonj}
\end{equation}
\smallskip By Lemma \ref{findim} there exists such an operator $T$ which
satisfies $\left\| T\right\| _{\vec{P}\rightarrow \vec{Y}}=c_{a}$.

\smallskip Evidently the functions $g_{1}$ and numbers $\beta _{1}$ and $%
\beta _{2}$ which are used in the formula defining $T$ must satisfy $%
|g_{1}(\xi )|\le c_{a}w_{0}(\xi )$ for a.e.\ $\xi \in [0,\alpha _{1})$ and $%
\sqrt{\beta _{1}^{2}+\beta _{2}^{2}}\le c_{a}\alpha _{2}$. Consequently,
substituting from (\ref{blonj}), we have
\begin{eqnarray*}
\alpha _{1} &=&\int_{[0,\alpha _{1})}g_{1}(\xi )d\xi +\beta _{1}\le
\int_{[0,\alpha _{1})}c_{a}w_{0}(\xi )d\xi +\sqrt{c_{a}^{2}\alpha
_{2}^{2}-\beta _{2}^{2}}= \\
&=&c_{a}(\sqrt{\alpha _{1}^{2}+\alpha _{2}^{2}}-\alpha _{2})+\alpha _{2}%
\sqrt{c_{a}^{2}-1}=c_{a}(1-\alpha _{2})+\alpha _{2}\sqrt{c_{a}^{2}-1}.
\end{eqnarray*}
In the special case where $a=\pi /6$, i.e.\ when $\alpha _{1}=\sqrt{3}/2$
and $\alpha _{2}=1/2$, the previous inequalities immediately imply that
\begin{equation*}
\sqrt{3}\le c_{\pi /6}+\sqrt{c_{\pi /6}^{2}-1}.
\end{equation*}
This is false if $c_{\pi /6}<2/\sqrt{3}$. I.e., we have shown that
\begin{equation}
c_{\pi /6}\ge 2/\sqrt{3}.  \label{zzz}
\end{equation}

We shall now prove that $c_{a}\leq 2/\sqrt{3}$ for all $a\in (0,\pi /2)$.
Having chosen such a value of $a$, we set $\alpha _{1}=\cos a$ and $\alpha
_{2}=\sin a$. Since
\begin{equation*}
\alpha _{1}+\alpha _{2}=\sqrt{\alpha _{1}^{2}+2\alpha _{1}\alpha _{2}+\alpha
_{2}^{2}}\mathbf{>}\sqrt{\alpha _{1}^{2}+\alpha _{2}^{2}}=1,
\end{equation*}
we have that $\frac{\alpha _{1}}{1-\alpha _{2}}>1$. It is clear that the
function $\phi (x):=\frac{\alpha _{1}-x}{1-\alpha _{2}}$ decreases from $%
\frac{\alpha _{1}}{1-\alpha _{2}}$ to $1$ on the interval $I=[0,\alpha
_{1}+\alpha _{2}-1]$. This in turn means that the continuous function $\psi
(x):=\alpha _{2}\sqrt{\phi ^{2}(x)-1}-x$ is also decreasing on the same
interval. Since $\psi (0)>0$ and $\psi (\alpha _{1}+\alpha _{2}-1)=1-\alpha
_{1}-\alpha _{2}<0$, there exists a number $\beta _{1}$ in the interior of $%
I $ such that $\psi (\beta _{1})=0$, i.e.\

\begin{equation}
\beta _{1}=\alpha _{2}\sqrt{\biggl(\frac{\alpha _{1}-\beta _{1}}{1-\alpha
_{2}}\biggr)^{2}-1}.  \label{propb}
\end{equation}
We shall use this number in the formula (\ref{deft}) to define an operator $%
T:\vec{P}\rightarrow \vec{Y}$ where we choose the other numbers and
functions in the formula by setting $g_{2}(\xi )=0$ (as we are obliged to
do) and also
\begin{equation}\label{g1def}
g_{1}(\xi )=\frac{\alpha _{1}-\beta _{1}}{1-\alpha _{2}}w_{0}(\xi )\text{
for all }\xi \in [0,\alpha _{1})\text{ and }\beta _{2}=\alpha _{2}.
\end{equation}

Observe that, with these definitions,
\begin{equation*}
T\chi _{[0,\alpha _{1}]}=\left( \frac{\alpha _{1}-\beta _{1}}{1-\alpha _{2}}%
(1-\alpha _{2})+\beta _{1},\alpha _{2}\right) =(\alpha _{1},\alpha _{2}),
\end{equation*}
i.e.\ the quantity $C_{a}:=\Vert T\Vert _{\vec{P}\to \vec{Y}}$ belongs to $%
E_{a}^*$. In particular, $c_{a}\leq C_{a}$. But, in view of (\ref{ching}) and (%
\ref{propb}), we have
\begin{equation*}
C_{a}=\frac{\alpha _{1}-\beta _{1}}{1-\alpha _{2}}.
\end{equation*}
This in turn can be substituted in (\ref{propb}) to give
\begin{equation*}
\beta _{1}=\alpha _{2}\sqrt{C_{a}^{2}-1}
\end{equation*}
and so
\begin{equation*}
C_{a}=\frac{\alpha _{1}-\beta _{1}}{1-\alpha _{2}}=\frac{\alpha _{1}-\alpha
_{2}\sqrt{C_{a}^{2}-1}}{1-\alpha _{2}}.
\end{equation*}
We deduce that
\begin{equation}
C_{a}+\frac{\alpha _{2}}{1-\alpha _{2}}\sqrt{C_{a}^{2}-1}=\frac{\alpha _{1}}{%
1-\alpha _{2}}.  \label{ojoj}
\end{equation}
We claim that (\ref{ojoj}) implies that
\begin{equation}
C_{a}\le 2/\sqrt{3}.  \label{jo}
\end{equation}
If this is false, then
\begin{equation*}
\frac{\alpha _{1}}{1-\alpha _{2}}>\frac{2}{\sqrt{3}}+\frac{\alpha _{2}}{%
1-\alpha _{2}}\sqrt{\frac{4}{3}-1}=\frac{1}{\sqrt{3}}\left( 2+\frac{\alpha
_{2}}{1-\alpha _{2}}\right)
\end{equation*}
and so $\sqrt{3}\alpha _{1}>2(1-\alpha _{2})+\alpha _{2}=2-\alpha _{2}$.
Consequently, $3\alpha _{1}^{2}>4-4\alpha _{2}+\alpha _{2}^{2}$. Since $%
\alpha _{1}^{2}+\alpha _{2}^{2}=1$ it follows that $3-3\alpha
_{2}^{2}>4-4\alpha _{2}+\alpha _{2}^{2}$, i.e.\ that $4\alpha
_{2}^{2}-4\alpha _{2}+1<0$. But this cannot hold for any real number $\alpha
_{2}$. This contradiction establishes (\ref{jo}).

We immediately deduce that $c_{a}\le 2/\sqrt{3}$ for all $a\in (0,\pi /2)$.
Combining this with (\ref{zzz}) and (\ref{opin}) gives (\ref{cy}) and
completes the proof of the theorem. \qed

\subsection{\label{genr}Generalizations and further remarks.}

\smallskip

We have the following generalization of Theorem \ref{cfy}.

\begin{theorem}
Let $U$ and $V$ be nontrivial Hilbert spaces and consider the
couple $\vec{W}=(U\oplus V,U)$. Then $\gamma(\vec{W})=2/\sqrt{3}$.
\end{theorem}

The proof is very similar to the case of $\vec{Y}$. We sketch the
changes necessary to make the proof work in the general case.

Choose unit vectors $u\in U$ and $v\in V$ and an element
$\alpha=(\alpha_1,\alpha_2)\in {\mathbb R}^2$ such that
$\alpha_1^2+\alpha_2^2=1$. It is then easy to see that
$$K(t,\alpha_1u+\alpha_2v;\vec{W})=K(t,\alpha;\vec{Y})=K(t,f;\vec{P}),\quad t>0,$$
with $f=\chi_{[0,\alpha_1]}$ and $\vec{P}=(L^1_{w_0},L^1_{w_1})$
defined as before. For $w\in W_0+W_1$ let $c_{w}=c_w(\vec{W})$ be
the quantity defined by (\ref{rbb}). It follows from Remark
\ref{4lat} that
$$\gamma(\vec{W})=\sup\{c_{\alpha_1u+\alpha_2v}\},$$
the supremum being taken over all points $(\alpha_1,\alpha_2)$ of
the unit circle and all unit vectors $u\in U$, $v\in V$.

For fixed $u, v$ and $\alpha$ as above we now choose
$T:\vec{P}\to\vec{W}$ as
\begin{equation}\label{canrep}
Th=\biggl(\int_{[0,\alpha_1)}g_1(\xi)h(\xi)d\xi+\beta_1h(\alpha_1)\biggr)u+
\alpha_2h(\alpha_1)v\end{equation} where the functions $g_1$ and
the number $\beta_1$ are defined by (\ref{propb}) and
(\ref{g1def}). Clearly, $Tf=\alpha_1u+\alpha_2v$. Moreover, as in
the case for $\vec{Y}$, one verifies that this operator $T$
satisfies
$$\left\| T\right\| _{\vec{P}\rightarrow \vec{W}}=\max \left\{ \mathop{\mathrm{ess%
\ sup}}_{\xi \in [0,\alpha _{1})}\left| \frac{g_{1}(\xi )}{w_{0}(\xi )}%
\right| ,\frac{\sqrt{\beta _{1}^{2}+\alpha_{2}^{2}}}{\alpha
_{2}}\right\}.$$ By the reasoning at the end of the proof of
Theorem \ref{cfy}, we now obtain that $c_{\alpha_1u+\alpha_2v}\le
\|T\|_{\vec{P}\to \vec{W}}\le 2/\sqrt{3}$, proving that
$\gamma(\vec{W})\le 2/\sqrt{3}$.

In order to prove the reverse inequality, we observe that an
arbitrary operator $S:\vec{P}\to \vec{W}$ such that
$Sf=\alpha_1u+\alpha_2v$ can be represented in the form
$$Sh=\biggl(\int_{[0,\alpha_1)}G_1(\xi)h(\xi)d\xi+B_1h(\alpha_1)\biggr)\oplus
(\alpha_2h(\alpha_1)v)$$ where $h\in P_0+P_1$ and $G_1\in
L^\infty([0,\alpha_1), U)$ and $B_1\in U$. Putting
$g_1(\xi)=(G_1(\xi),u)$ and $\beta_1=(B_1,u)$, we obtain a
corresponding operator $T$ of the form (\ref{canrep}) also
satisfying $Tf=\alpha_1u+\alpha_2v$ and such that
$\|T\|_{\vec{P}\to\vec{W}}\le\|S\|_{\vec{P}\to\vec{W}}$. Now, in
the case when $\alpha_1=\sqrt{3}/2$ and $\alpha_2=1/2$, the
estimate $\|T\|_{\vec{P}\to\vec{W}}\ge 2/\sqrt{3}$ follows exactly
as in the case for $\vec{Y}$. \hfill $\qed$

\smallskip

It seems plausible that couples of the above form are extremal
amongst all Hilbert couples in the sense that their
$K$-divisibility constant is maximal. Thus we have the following
open question.

\begin{question} Is $\gamma(\vec{H})\le 2/\sqrt{3}$ for every Hilbert
couple $\vec{H}$?
\end{question}

\smallskip
For a comment related to this question, see Remark \ref{lcoup} below.

\smallskip

We now turn to some generalizations of our result in other
directions. These will include the following result:

\begin{theorem}\label{sprop}
Let $\vec{X}=(X_0,X_1)$ be a Banach couple such that $X_0$ is
two-dimensional and $X_1$ is a one-dimensional subspace of $X_0$.
Then $\gamma(\vec{X})\le 2\sqrt{2/3}$.
\end{theorem}

\begin{remark} Note that Shvartsman's couple
$\vec{S}$ \cite{shv}, where $S_0$ is ${\mathbb R}^2$ equipped with
the $\ell_\infty$-norm and $S_1$ the one-dimensional subspace of
${\mathbb R}^2$ whose unit ball makes an angle of $\pi/8$ with the
positive $x$-axis, is of the form occurring in Theorem \ref{sprop}
and satisfies $\gamma(\vec{S})\ge \frac {3+2\sqrt{2}}
{1+2\sqrt{2}}$.
\end{remark}

In order to prove Theorem \ref{sprop} we first need to introduce
some terminology and obtain some preliminary results.

\begin{definition}
\label{rigid}Let $\vec{A}=(A_{0},A_{1})$ and $\vec{B}=(B_{0},B_{1})$
be two Banach couples. A linear operator $T:A_0+A_1\to B_0+B_1$
which, for $j=0$ and $j=1$, is a one to one map of $A_j$ onto $B_j$
and satisfies $\left\| Ta\right\| _{B_{j}}=c_{j}\left\| a\right\|
_{A_{j}}$ for all $a\in A_j$ and some positive constant $c_j$, is
called a \textit{\textbf{rigid map}} of $\vec{A}$ onto $\vec{B}$. If
such a map exists, then we say that $\vec{B}$ is a
\textit{\textbf{rigid image}} of $\vec{A}$. (This is of course the
same as saying that $\vec{A}$ is a rigid image of $\vec{B}$.)
\end{definition}

\smallskip A classical and much used example of two couples which are rigid
images of each other, goes back to the paper \cite{sw} of Stein
and Weiss, where it was pointed out that, in the terminology of
Definition \ref{rigid}, any couple of weighted $L^{p}$ spaces
$\vec{B}=\left( L_{w_{0}}^{p_{0}}(\Omega ,\Sigma ,\mu
),L_{w_{1}}^{p_{1}}(\Omega ,\Sigma ,\mu )\right)$ on some measure
space $(\Omega ,\Sigma ,\mu )$ where $1\le
p_{0}<p_{1}\le \infty $, is a rigid image of an unweighted couple $\vec{A}%
=\left( L^{p_{0}}(\Omega ,\Sigma ,\nu ),L^{p_{1}}(\Omega ,\Sigma
,\nu )\right) $ for some other measure $\nu $ on the same measure
space.

\smallskip

\begin{fact}
\label{rkd}If $\vec{B}$ is a rigid image of $\vec{A}$ then $\gamma (\vec{B}%
)=\gamma (\vec{A})$. Furthermore we have that $c_a(\vec{A})=c_{ Ta
}(\vec{B})$ for all $a\in A_0+A_1$, where $T$ is a rigid map of
$\vec{A}$ onto $\vec{B}$.
\end{fact}

In order to prove Fact \ref{rkd}, we first note that standard
arguments show immediately that $K(t,Ta;\vec{B})= c_{0}K\left(
\frac{c_{0}t}{c_{1}},a;\vec{A}\right)$ for all $t>0$ and all $a\in
A_{0}+A_{1}$.

Put $b=Ta$ and suppose that $K(t,b;\vec{%
B})\le \sum_{n=1}^{\infty }\psi _{n}(t)$ for all $t>0$, where the functions $%
\psi _{n}:(0,\infty )\rightarrow (0,\infty )$ are all concave and $%
\sum_{n=1}^{\infty }\psi _{n}(1)<\infty $. Then $K(t,b;\vec{B})=c_{0}K\left( \frac{%
c_{0}t}{c_{1}},a;\vec{A}\right) \le \sum_{n=1}^{\infty }c_{0}\psi
_{n}\left( \frac{c_{0}t}{c_{1}}\right) $. Since $\phi
_{n}(t):=c_{0}^{-1}\psi
_{n}\left( \frac{c_{1}t}{c_{0}}\right)$ is concave for each $n$ and $%
\sum_{n=1}^{\infty }\phi _{n}(1)<\infty $ it follows from Theorem \ref{ec} that, for each $%
\epsilon >0$, there exists a sequence of elements $\left\{
a_{n}\right\} _{n\in \Bbb{N}}$ in $A_{0}+A_{1}$ such that
$a=\sum_{n=1}^{\infty }a_{n}$ with convergence in $A_{0}+A_{1}$ norm
and $K(t,a_{n};\vec{A})\le \left( c_a(\vec{A}) +\epsilon \right)
\phi _{n}(t)$ for all $t>0$ and all $n\in
\Bbb{N}$. If we set $b_{n}=Ta_{n}$ for each $n$ then it is clear that $%
b=\sum_{n=1}^{\infty }b_{n}$ with convergence in $B_{0}+B_{1}$ norm
and
\begin{equation*}
K(t,b_{n};\vec{A})=c_{0}K\left(
\frac{c_{0}t}{c_{1}},a_{n};\vec{A}\right)
\le \left( c_a(\vec{A}) +\epsilon \right) c_{0}\phi _{n}\left( \frac{%
c_{0}t}{c_{1}}\right) =\left( c_a(\vec{A}) +\epsilon \right) \psi
_{n}\left( t\right)
\end{equation*}
for all $t>0$ and all $n\in \Bbb{N}$. This shows that
$c_{b}(\vec{B}) \le
c_a(\vec{A}) +\epsilon $ for each positive $\epsilon $. It follows that $%
c_{b}(\vec{B}) \le c_a(\vec{A})$ and of course an analogous argument
using $T^{-1}$ in place of $T$ shows that $c_a(\vec{A}) \le c_{b}(\vec{B}%
)$. This finishes the proof of Fact \ref{rkd}.

\smallskip

Suppose that $\vec{X}=(X_0,X_1)$ satisfies the hypotheses of Theorem
\ref{sprop} and that, furthermore, $X_0$ is also a Hilbert space.
Then it is easy to see that $\vec{X}$ is a rigid image of
$\vec{Y}=(\ell^2_2,\ell^2_1)$, and consequently $\gamma(\vec{X}) =2/\sqrt{3}$.
(More explicitly, suppose that $\{e_1,e_2\}$ is an orthonormal
basis for $X_0$. Then, for some constants $\alpha$ and $\beta$ we have
$X_1=\{\alpha t e_1+\beta t e_2:t\in \Bbb{R} \}$ and
$\| \alpha e_1+\beta e_2 \|_{X_1}=1$. Now let $V_0=\Bbb{R}^2$
with $\|(x,y)\|_{V_0}=\|x e_1+ye_2\|_{X_0}=\sqrt{x^2+y^2}$ and
let $V_1=\{(\alpha t,\beta t):t\in \Bbb{R}\}$ with
$\| (\alpha t,\beta t)\|_{V_1}=\| \alpha t e_1+\beta e_2 t \|_{X_1}=|t|$.
The linear map $T:\vec{V}\to \vec{X}$ defined by
$T(x,y)=xe_1+ye_2$ shows that $\vec{V}$ and $\vec{X}$ are rigid images
of each other. Then a suitable map of rotation in $\Bbb{R}^2$ which moves
the point $(\alpha ,\beta )$ to $(\sqrt{\alpha^2+\beta^2},0)$
shows that $\vec{V}$ is a rigid image of $\vec{Y}$.)

\smallskip

The classical Banach-Mazur distance between Banach space has a
counterpart for Banach couples. We have the following definition.

\begin{definition}\label{banachmazur}
Let $\vec{A}=(A_0,A_1)$ and $\vec{B}=(B_0,B_1)$ be Banach couples.
If $A_j$ is isomorphic to $B_j$ for $j=0,1$, then the
\textit{\textbf{Banach-Mazur distance}} between $\vec{A}$ and
$\vec{B}$ is defined by
$$d(\vec{A};\vec{B})=\inf\{\|T\|_{\vec{A}\to\vec{B}}\|T^{-1}\|_{\vec{B}\to\vec{A}}\},$$
the infimum being taken over all linear isomorphisms
$T:\vec{A}\to\vec{B}$. Otherwise $d(\vec{A};\vec{B})=\infty$.
\end{definition}

Several of the general results which now follow should probably be
considered as belonging to the folklore of interpolation theory.
For example, they should be compared with Section 3 of Brudnyi and
Shteinberg \cite{bs}, where similar notions and results are
discussed.

\begin{proposition} \label{bcont}
Let $\vec{A}$ and $\vec{B}$ be non-zero Banach couples. Then
\begin{equation}\label{gcont}
\gamma(\vec{A})\le\gamma(\vec{B})d(\vec{A};\vec{B}).
\end{equation}
In particular, $\gamma$ is a bounded continuous function on the
category of Banach couples endowed with the Banach-Mazur metric.
\end{proposition}

Before we prove this proposition let us show how Theorem
\ref{sprop} follows from it and Theorem \ref{cfy}:

Use John's theorem to choose a two-dimensional Hilbert space $Z_0$
such that
$\|\cdot\|_{Z_0}\le\|\cdot\|_{X_0}\le\sqrt{2}\|\cdot\|_{Z_0}$ and
let $Z_1=X_1$. The couple $\vec{Z}$ is then a Hilbert couple such
that $d(\vec{S};\vec{Z})\le\sqrt{2}$. As explained above,
$\vec{Z}$ is isometric to a rigid image of the couple $\vec{Y}$
and so we have $\gamma(\vec{Z})=2/\sqrt{3}$ which proves Theorem
\ref{sprop}.

It remains to prove Proposition \ref{bcont}. The fact that
$\gamma$ is bounded is of course the Brudnyi-Krugljak theorem
(Theorem \ref{bkthm}), so we will only need to prove
(\ref{gcont}).

In fact, we will deduce (\ref{gcont}) from a more general
proposition. We will first require yet another definition: (Cf.\
\cite{bs}.)

\begin{definition}\label{bscc}
Let ${\mathbb K}$ be either ${\mathbb R}$ or ${\mathbb C}$ and
assume in the following that all Banach spaces are over the field
${\mathbb K}$.

Let $C$ be a non-negative constant. Two couples $\vec{A}$,
$\vec{B}$ are \textit{\textbf{relative $C$-monotonic couples}} if
for every $\eps>0$, all $\alpha\in A_0+A_1$ and $\beta\in B_0+B_1$
such that
\begin{equation}\label{kkond}
K(t,\beta;\vec{B})\le K(t,\alpha;\vec{A}),\quad t>0
\end{equation}
there exists a ${\mathbb K}$-linear operator
$T=T_\eps:\vec{A}\to\vec{B}$ such that
$$T\alpha=\beta\text{ and
}\|T\|_{\vec{A}\to\vec{B}}< C+\eps.$$
The smallest constant $C$
satisfying this implication is called the
\textit{\textbf{Calder\'{o}n constant}} relative to $\vec{A}$ and
$\vec{B}$ and is denoted by $c(\vec{A};\vec{B})$. We also put
$$c_n({\mathbb K})=\sup\{c(\vec{A};\vec{B}):\dim_{\mathbb K}(A_i)\le n\text{ and
}\dim_{\mathbb K}(B_i)\le n,\, i=0,1\}.$$
\end{definition}

\begin{proposition}\label{cccont} Let $\vec{A}^i$ and $\vec{B}^i$ be non-zero Banach
couples for $i=1,2$.
Then
$$c(\vec{A}^1;\vec{B}^1)\le
d(\vec{A}^1;\vec{A}^2)c(\vec{A}^2;\vec{B}^2)d(\vec{B}^1;\vec{B}^2).$$
\end{proposition}

\begin{proof} We may assume that both of the Banach-Mazur distances above
are finite, because otherwise the statement is trivial.

Take $\alpha\in A^1_0+A_1^1$ and $\beta\in B_0^1+B_1^1$ such that
$K(t,\beta;\vec{B}^1)\le K(t,\alpha;\vec{A}^1)$ for all $t>0$. Let
$T_A:\vec{A}^1\to \vec{A}^2$ and $T_B:\vec{B}^1\to\vec{B}^2$ be
isomorphisms such that
$\|T_A\|\|T_A^{-1}\|<d(\vec{A}^1;\vec{A}^2)+\epsilon$ and
$\|T_B\|\|T_B^{-1}\|<d(\vec{B}^1;\vec{B}^2)+\epsilon$. It follows
that
$$\|T_B\|^{-1}K(t,T_B(\beta);\vec{B}^2)\le K(t,\alpha;\vec{A}^1)\le
\|T_A^{-1}\|K(t,T_A(\alpha);\vec{A}^2)$$ for all $t>0$. Take
$\epsilon>0$. It then follows that there exists an operator
$T_0:\vec{A}^2\to\vec{B}^2$ such that
$T_0(T_A(\alpha))=T_B(\beta)$ of norm at most
$(c(\vec{A}^2;\vec{B}^2)+\epsilon)\|T_B\|\|T_A^{-1}\|$. The
operator $T:\vec{A}^1\to\vec{B}^1$ defined by $T=T_B^{-1}\circ
T_0\circ T_A$ then fulfills $T(\alpha)=\beta$ and $\|T\|\le
\|T_B^{-1}\|(c(\vec{A}^2;\vec{B}^2)+\epsilon)\|T_B\|\|T_A^{-1}\|\|T_A\|<d(\vec{A}^1;\vec{A}^2)
c(\vec{A}^2;\vec{B}^2)d(\vec{B}^1;\vec{B}^2)+O(\epsilon)$.
\end{proof}

\begin{proof}[Proof of Proposition \ref{bcont}]
Fix a Banach couple $\vec{A}$. By Theorem \ref{ec} we have
$$\gamma(\vec{A})=\sup \{c(\vec{P};\vec{A})\}$$
over weighted $L^1$-couples $\vec{P}$. But Proposition \ref{cccont}
yields that for each particular weighted $L^1$-couple $\vec{P}$ we
have
$$
c(\vec{P};\vec{A})\le c(\vec{P};\vec{B})d(\vec{A};\vec{B}).
$$
The inequality (\ref{gcont}) follows by taking the supremum over
all weighted $L^1$ couples $\vec{P}$.
\end{proof}

\begin{remark} \label{lcoup}
Let $\vec{H}$ be a finite-dimensional Hilbert couple. Then it is
easy to see that there exists a finite sequence
$\lambda=(\lambda_i)_{i=1}^n\subset [0,\infty]$ such that
$\vec{H}$ is isometric to the weighted $\ell^2$-couple
$(\ell_n^2,\ell_n^2(\lambda))$. A generalization of this statement
to the case of infinite-dimensional Hilbert couples has been given
by Sedaev \cite{sed}. By this latter observation, the
interpolation of Hilbert couples becomes essentially the same as
that of weighted $\ell^2$-couples. (Cf. also \cite{do} and
\cite{ameur}, \cite{a2}.)
\end{remark}

\subsection{\label{caldconst}Calder\'on constants
for finite dimensional couples.}

Since we have had to introduce and use relative Calder\'on
constants in the previous subsection, it is now convenient for us
to make a slight digression and prove the following theorem
estimating the size of relative Calder\'{o}n constants for couples
of a given finite dimension. This result is closely related to
Theorem 3.1 of \cite{bs}. The method of proof is is very similar
to that of \cite{bs} Section 3.

\begin{theorem} \label{prexis} $c_n({\mathbb C})=n$  and $n/\sqrt{2}\le
c_n({\mathbb R})\le n$ for all $n\in{\mathbb N}$.
\end{theorem}

\begin{remark} In \cite{bs}, Brudnyi and Shteinberg
introduce the quantity $\varkappa _n$ defined by
$$\varkappa _n=\sup\{c(\vec{A};\vec{A}):\dim(A_i)\le n \text{ for }
i=0,1\},$$ where the supremum is taken with respect to Banach
couples over the reals. In Theorem 3.1 of \cite{bs} they show that
$n/2\sqrt{2}\le \varkappa _n\le n\sqrt{2}$. Since of course
$\varkappa _n \le c_n(\Bbb{R})$, our result provides a somewhat
better upper estimate for $\varkappa _n$.
\end{remark}

\begin{proof}[Proof of Theorem \ref{prexis}]

`` $\le$": Let $\vec{A}$ and $\vec{B}$ be couples such that all
the spaces $A_i$ and $B_i$ are of dimension at most $n$ (scalars
can be real or complex). Let $\alpha\in A_0+A_1$ and $\beta\in
B_0+B_1$ be elements satisfying (\ref{kkond}). Use John's theorem
to find Hilbert spaces $H_i$ and $K_i$ such that
$d(\vec{A};\vec{H})\le\sqrt{n}$ and
$d(\vec{B};\vec{K})\le\sqrt{n}$. By Proposition \ref{cccont}
$$c(\vec{A};\vec{B})\le nc(\vec{H};\vec{K})$$
But Hilbert couples are exact relative Calder\'{o}n couples, i.e.,
$c(\vec{H};\vec{K})\le 1$ by Theorem 2.2 of \cite{ameur}. Thus
$c(\vec{A};\vec{B})\le n$.

`` $\ge$": This is a straightforward adaptation of the elegant
arguments given in \cite{bs}, Section 3.

First assume complex scalars and define the space
$\ell^{p,r}_n(q)$ for suitable fixed values of $p, q$ and $r$ by
the norm
$$\|x\|_{\ell^{p,r}_n(q)}^p=\sum_{k=1}^n |q^{-kr}x_k|^p,\quad x=(x_k)_1^n\in
{\mathbb C}^n.$$ For fixed $p$ and $q$ we also define the couple
$\vec{\ell}^p_n(q)=(\ell^{p,0}_n(q),\ell^{p,1}_n(q))$. (The usual
conventions apply for the case $p=\infty$.)

Choose a fixed $q>1$ and put $h=(\sqrt{q},\sqrt{q}^2,
\ldots,\sqrt{q}^{n})\in {\mathbb C}^n$. As is shown in \cite{bs},
we have
\begin{equation}\label{jinxx}
K(t,h;\vec{\ell}^1_n(q))=\sum_{k=1}^n q^{k/2}\min\{1,q^{-k}t\}\le
\frac {\sqrt{q}-1} {\sqrt{q}+1}K(t,h;\vec{\ell}_n^\infty(q)).
\end{equation}
(It is convenient to first prove the inequality in the cases $t=q^i$,
and then use the concavity of the $K$-functional.)

By (\ref{jinxx}) there exists an operator
$T:\vec{\ell}^\infty_n(q)\to \vec{\ell}^1_n(q)$ such that $T(h)=h$
and $\|T\|\le \frac {\sqrt{q}+1}
{\sqrt{q}-1}c(\vec{\ell}^\infty_n(q);\vec{\ell}^1_n(q))$.

Since we are assuming complex scalars, the Riesz--Thorin theorem
can be applied. It yields that
$$\|T\|_{\ell^{\infty,1/2}_n(q)\to\ell^{1,1/2}_n(q)}\le \frac
{\sqrt{q}+1}
{\sqrt{q}-1}c(\vec{\ell}^\infty_n(q);\vec{\ell}^1_n(q)).$$ This in
turn yields
$$n=\|h\|_{\ell^{1,1/2}_n(q)}\le\frac
{\sqrt{q}+1}
{\sqrt{q}-1}c(\vec{\ell}^\infty_n(q);\vec{\ell}^1_n(q))\|h\|_{\ell^{\infty,1/2}_n(q)}=\frac
{\sqrt{q}+1}
{\sqrt{q}-1}c(\vec{\ell}^\infty_n(q);\vec{\ell}^1_n(q)).$$ It
follows that $c_n({\mathbb C})\ge
c(\vec{\ell}_n^\infty(q);\vec{\ell}^1_n(q))\ge n\frac {\sqrt{q}-1}
{\sqrt{q}+1}$. Since $q$ can be chosen arbitrarily large, this
gives $c_n({\mathbb C})\ge n$. The modifications necessary to
treat the real case are carried out as in \cite{bs}.
\end{proof}

We end this subsection with an open question.

\begin{question} Is $c_n({\mathbb R})=n$?
\end{question}

\smallskip

\subsection{\label{hcmore}On the case of a regular two dimensional
Hilbert couple.}

Let $r$ be a positive number and let $\vec{G}=(G_{0},G_{1})$ be
the couple for which $G_{0}=\ell _{2}^{2}$ and $G_{1}$ is the
weighted version of $\ell _{2}^{2}$ with norm $\left\|
(x,y)\right\| _{G_{1}}=\sqrt{x^{2}+ry^{2}}$.

\smallskip In this subsection we will prove a rather simple estimate: $\gamma (\vec{G})<\sqrt{2}$.

Let us remark first that in the trivial case where $r=1$ we obtain
$\gamma(\ell_2^2,\ell_2^2)=1$. In the general case, Proposition
\ref{bcont} yields that $\gamma(\vec{G})$ is a continuous function
of $r$ and $\gamma(\vec{G})\le\max(\sqrt{r},1/\sqrt{r})$. (This is
because the Banach--Mazur distance between $\vec{G}$ and
$(\ell_2^2,\ell_2^2)$ is $\max(\sqrt{r},1/\sqrt{r})$.)

Fix a point $\alpha =(b,c)=(\cos a,\sin a)\in G_{0}+G_{1}$ where
$a\in [0,2\pi )$. In fact, by Remark \ref{4lat}, we only need to
consider the case where $a\in [0,\pi /2]$.

\smallskip

We will look for a parametric representation of the curve which is
the boundary $\partial \Gamma (\alpha )$ of the Gagliardo diagram
of $\alpha $.

First let us fix some $t>0$ and determine the point $z=(x,y)$ for
which the infimum $K_{2}(t,\alpha ;G_{0},G_{1})^2=\inf_{z\in
\Bbb{R}^{2}}\left\| z\right\| _{G_{0}}^{2}+t\left\| \alpha
-z\right\| _{G_{1}}^{2}$ is attained. The point which we are looking
for is of course the unique critical point of the function $\phi
(x,y)=x^{2}+y^{2}+t\left( x-b\right) ^{2}+tr(y-c)^{2}$, i.e.\ $x=\frac{tb}{%
1+t}$ and $y=\frac{trc}{1+tr}$.

It is clear that, for this choice of $z$, the point $\left(
\left\| z\right\| _{G_{0}}^{2},\left\| \alpha -z\right\|
_{G_{1}}^{2}\right) $\smallskip belongs to $\partial \Gamma
(\alpha )$, and
that, furthermore, as $t$ ranges over $(0,\infty )$ we obtain all points of $%
\partial \Gamma (\alpha )\cap \left\{ (x_{0},x_{1}):x_{0}>0,x_{1}>0\right\} $
in this way. We note that \smallskip
$b-x=\frac{b+tb-tb}{1+t}=\frac{b}{1+t}$ and
$c-y=\frac{c+trc-trc}{1+tr}=\frac{c}{1+tr}$. It follows that
\begin{equation}
\partial \Gamma (\alpha )\cap \left\{ (x_{0},x_{1}):x_{0}>0,x_{1}>0\right\}
=\left\{ \left( \gamma _{0}(t),\gamma _{1}(t)\right) :0<t<\infty
\right\} . \label{jonqz}
\end{equation}
\smallskip where the functions $\gamma _{0}$ and $\gamma _{1}$ are given by
\begin{equation*}  
\gamma _{0}(t)=t\sqrt{\frac{b^{2}}{(1+t)^{2}}+\frac{r^{2}c^{2}}{(1+tr)^{2}}}%
\text{ and }\gamma _{1}(t)=\sqrt{\frac{b^{2}}{(1+t)^{2}}+\frac{rc^{2}}{%
\left( 1+tr\right) ^{2}}}\text{ for all }t\in (0,\infty ).
\end{equation*}
Obviously $\gamma _{1}(t)$ is a strictly decreasing function of $t$. Since $%
\gamma_{0}(1/t)^2=\frac{b^{2}}{(t+1)^{2}}+\frac{r^{2}c^{2}}{(t+r)^{2}}$
it is also clear that $\gamma _{0}(t)$ is a strictly increasing function of $%
t$.

Considering the limits of $\gamma _{0}$ and $\gamma _{1}$ as $t$
tends to $0$ and to $\infty $, we deduce that $\partial \Gamma
(\alpha )$ is the union of the curve specified in (\ref{jonqz})
with the two rays on the coordinate axes
\begin{equation}\label{glimits}
\left\{ (0,v):\sqrt{b^{2}+rc^{2}}\le v<\infty \right\} \text{ and
}\left\{ (v,0):1\le v<\infty \right\} .
\end{equation}
Next we define two functions $w_{0}$ and $w_{1}$ by
$w_{0}(t):=\gamma _{0}^{\prime }(t)$ and $w_{1}(t):=-\gamma
_{1}^{\prime }(t)$ for all $t\in (0,\infty )$. These will turn out
to be convenient weight functions to use in a couple of weighted
$L^{1}$ spaces on $(0,\infty )$ as an essential step for
calculating $\gamma ( \vec{G}) $. We note that (\ref{glimits})
implies
\begin{equation}\label{wint}
\int_0^\infty w_0(t)dt=1\quad {\rm and}\quad \int_0^\infty
w_1(t)dt=\sqrt{b^2+rc^2}.
\end{equation}

We will see that routine calculations
show that $w_0$ and $w_1$ are given explicitly by
\begin{equation}
w_{j}(t)=\frac{\displaystyle \frac{b^{2}}{(1+t)^{3}}+\frac{r^{2}c^{2}}{%
(1+rt)^{3}}}{\sqrt{\displaystyle \frac{b^{2}}{(1+t)^{2}}+\frac{r^{2-j}c^{2}}{%
\left( 1+rt\right) ^{2}}}}\text{ for }j=0,1\text{ and }t\in
(0,\infty ). \label{lebc}
\end{equation}
The proof of this in the case $j=1$ is immediate. For the case
$j=0$ we can first observe that
\begin{equation*}
\gamma _{0}^{\prime }(1/t)\cdot
\frac{1}{t^{2}}=-\frac{d}{dt}\left( \gamma
_{0}(1/t)\right) =\frac{\frac{b^{2}}{(1+t)^{3}}+\frac{r^{2}c^{2}}{(t+r)^{3}}%
}{\sqrt{\frac{b^{2}}{(1+t)^{2}}+\frac{r^{2}c^{2}}{\left(
t+r\right) ^{2}}}}
\end{equation*}
which implies that
\begin{equation*}
w_{0}(1/t)=\frac{\frac{t^{3}b^{2}}{(1+t)^{3}}+\frac{t^{3}r^{2}c^{2}}{%
(t+r)^{3}}}{t\sqrt{\frac{b^{2}}{(1+t)^{2}}+\frac{r^{2}c^{2}}{\left(
t+r\right) ^{2}}}}=\frac{\frac{b^{2}}{(1/t+1)^{3}}+\frac{r^{2}c^{2}}{%
(1+r/t)^{3}}}{\sqrt{\frac{b^{2}}{(1/t+1)^{2}}+\frac{r^{2}c^{2}}{\left(
1+r/t\right) ^{2}}}}
\end{equation*}
which immediately gives (\ref{lebc}) for $j=0$.

\smallskip Note that $w_{0}$ and $w_{1}$ are both strictly positive on $%
(0,\infty )$. \smallskip

We will use the couple $\vec{P}=(P_{0},P_{1})$ of weighted $L^{1}$
spaces
on the measure space $(0,\infty )$ (equipped with Lebesgue measure) where $%
P_{0}=L_{w_{0}}^{1}$ and $P_{1}=L_{w_{1}}^{1}$. Let $f$ be the
function which equals $1$ identically on $(0,\infty )$. We will show
that
\begin{equation}\label{ccheck}
K(t,f;\vec{P})=K(t,\alpha ;\vec{G})\text{ for all }t>0.
\end{equation}

For each $t>0$ it is well known and very easy to check that
\begin{equation}
K(t,f;\vec{P})=\int_{0}^{\infty }\min \left\{
w_{0}(s),tw_{1}(s)\right\} ds \label{usuw}
\end{equation}
and that an optimal decomposition $f=f_{0,t}+f_{1,t}$, for which the
infimum in the calculation of (\ref{usuw}) is attained, is given by
$f_{0,t}=\chi _{E_{t}}$ and $f_{1,t}=\chi _{(0,\infty )\backslash
E_{t}}$, where
\begin{equation}\label{defte}E_{t}=\left\{ s>0:w_{0}(s)<tw_{1}(s)\right\}.
\end{equation}

We need to consider the function
\begin{equation}
\frac{w_{0}(t)^2}{w_{1}(t)^2} =\frac{\frac{b^{2}}{(1+t)^{2}}+\frac{rc^{2}%
}{\left( 1+tr\right)
^{2}}}{\frac{b^{2}}{(t+1)^{2}}+\frac{r^{2}c^{2}}{\left(
1+rt\right) ^{2}}}=1+\frac{(r-r^{2})c^{2}}{b^{2}\left(
r+\frac{1-r}{t+1}\right) ^{2}+r^{2}c^{2}}.  \label{skm}
\end{equation}
In the trivial cases where $(b,c)$ is either $(0,1)$ or $(0,1)$ this
is a constant function, and it a simple matter to check that
(\ref{ccheck}) holds. (In the first case the $K$-functionals on the
left and right sides of (\ref{ccheck}) both equal $\min\{1,t\}$ and
in the second case they both equal $\min\{1,t\sqrt{r}\}$).

\smallskip In the
remaining non-trivial case when $b$ and $c$ are both non zero it is
easy to see from (\ref{skm}) that, for any $r\in (0,\infty )$ with
$r\ne 1$,
\begin{equation}
\frac{w_{0}(t)}{w_{1}(t)}\text{ is a strictly increasing
continuous function of }t\text{ on }(0,\infty ).  \label{fim}
\end{equation}
(The two cases $r<1$ and $r>1$ have to be considered separately.)

We introduce and calculate two ``limiting'' values of $t$ by
setting
\begin{equation}
t_{0}^2:=\lim_{s\rightarrow 0}\frac{w_{0}(s)^2}{w_{1}(s)^2}=1+\frac{%
(r-r^{2})c^{2}}{b^{2}+r^{2}c^{2}}=\frac {b^2+rc^2} {b^2+r^2c^2}
\label{lol}
\end{equation}
and
\begin{equation}
t_{\infty }^2:=\lim_{s\rightarrow \infty }\frac{w_{0}(s)^2}{w_{1}(s)^2}=1+%
\frac{(r-r^{2})c^{2}}{r^{2}(b^{2}+c^{2})}=b^2+c^2/r. \label{upl}
\end{equation}

The property (\ref{fim}) implies that the set $E_t$ defined in
(\ref{defte}) is an open interval of the form $E_{t}=$ $(0,u(t))$,
where $u$ is a
non decreasing function of $t$. By (\ref{lol}) and (\ref{upl}%
) we see that $u(t)=0$ for $t\le t_{0}$ and $u(t)=\infty $ for
$t\ge t_{\infty }$, and, for each $t\in (t_{0},t_{\infty })$,
$u(t)$ is the unique number in $(0,\infty )$ for which
$w_{0}(u(t))/w_{1}(u(t))=t$.

We can now deduce that, for $t\in (t_{0},t_{\infty })$, $\left\|
f_{0,t}\right\|
_{P_{0}}=\int_{0}^{u(t)}w_{0}(s)ds=\int_{0}^{u(t)}\gamma
_{0}^{\prime }(s)ds=\gamma _{0}(u(t))-\gamma _{0}(0)=\gamma _{0}(u(t))$ and $%
\left\| f_{1,t}\right\| _{P_{1}}=\int_{u(t)}^{\infty
}w_{1}(s)ds=-\int_{u(t)}^{\infty }\gamma _{1}^{\prime
}(s)ds=\gamma _{1}(u(t))-\lim_{r\rightarrow \infty }\gamma
_{1}(r)=\gamma _{1}(u(t))$.

This shows that, as $t$ ranges over the interval $(t_{0},t_{\infty
})$, the point $\left( \left\| f_{0,t}\right\| _{P_{0}},\left\|
f_{1,t}\right\| _{P_{1}}\right)$ ranges over the curve
(\ref{jonqz}), i.e., $\Gamma(f)=\Gamma(\alpha)$. By the well-known
relation between $K$-functionals and Gagliardo diagrams, (see
\cite{bl}, sect. 7,1), this implies that (\ref{ccheck}) holds.

\smallskip It is clear that every bounded operator $T:\vec{P}\rightarrow
\vec{G}$ uniquely determines and is uniquely determined by a
suitable pair of (equivalence classes of) measurable functions
$g_{j}:(0,\infty )\rightarrow \Bbb{R}$ for $j=0,1$, via the
formula
\begin{equation}
Th=\left( \int_{0}^{\infty }g_{0}(s)h(s)ds,\int_{0}^{\infty
}g_{1}(s)h(s)ds\right) \text{ for all }h\in
L_{w_{0}}^{1}+L_{w_{1}}^{1}. \label{kewc}
\end{equation}
When it is necessary to explicitly indicate the connection between
the operator $T$ and the functions $g_{0}$ and $g_{1}$ which
define it via (\ref {kewc}), we will use the notation
$T_{g_{0},g_{1}}$ in place of $T$.

Of course we need to be more explicit about the conditions that the
functions $g_{0}$ and $g_{1}$ must satisfy. Straightforward
arguments (exactly like the proof below
of the equivalence of conditions (\ref{jse}) and (%
\ref{zjse})) using the Lebesgue differentiation theorem and a
suitable form of Minkowski's or Schwartz' inequality, show that
the norm of $T$ is given by
\begin{equation}
\left\| T\right\| _{\vec{P}\rightarrow \vec{G}}=\max_{j=0,1}\left\{ \mathrm{%
ess}\sup_{(0,\infty
)}\frac{\sqrt{g_{0}^{2}+r^{j}g_{1}^{2}}}{w_{j}}\right\}
\label{oinq}
\end{equation}
and so $g_{0}$ and $g_{1}$ must be such that this expression in
finite.

\begin{remark}\label{rwlog} For our purposes, we can without loss
of generality assume that $r>1$, since for each $r<1$, the couple
$\vec{G}$ is a rigid image of the corresponding couple where $r$ has
been replaced by $1/r$. (Use Fact \ref{rkd} and the rigid map
$(x,y)\mapsto (y/\sqrt{r},x/\sqrt{r})$.)
\end{remark}

Now we will consider the class $\mathcal{T}=\mathcal{T}_{b,c}$ of
all bounded operators $T:\vec{P}\rightarrow \vec{G}$ which satisfy
$Tf=\alpha $ and consider the quantity
$c_a=c_a(\vec{G})=\inf\{\|T\|:T\in \mathcal{T} \}$. We
first make
a simple observation:

\begin{proposition} \label{sq2prop} We have $c_0=c_{\pi/2}=1$ and if $a\in (0,\pi/2)$ then $c_a<\sqrt{1+b^2}$.
In particular, $\gamma(\vec{G})<\sqrt{2}$.
\end{proposition}

\begin{proof}
By Remark \ref{rwlog} we can and will assume that $r>1$.

If $a=0$, i.e., if $(b,c)=(1,0)$, then we have that
$w_0^2(t)=w_1^2(t)=\frac 1 {(1+t)^2}$ and the operator
$T=T_{g_0,g_1}$ defined by $g_0(s)=\frac 1 {(1+s)^2}$ and
$g_1(s)=0$ satisfies $Tf=(b,c)$ and $\|T\|=1$. Thus $c_0=1$. The
proof of the fact that $c_{\pi/2}=1$ is equally simple. It uses
the functions $g_0(s)=0$ and $g_1(s)=\frac r {(1+rs)^2}$.

Now let $a\in (0,\pi/2)$. We claim that it suffices to consider
the operator $T=T_{g_0,g_1}$ given by $g_0(s)=bw_0(s)$ and
$g_1(s)=cw_1(s)/\sqrt{b^2+rc^2}$. Indeed $Tf=\alpha$ by
(\ref{wint}), and furthermore, by (\ref{lol}),
$$\frac {g_0^2(s)+g_1^2(s)} {w_0^2(s)}=b^2+\frac {c^2}
{b^2+rc^2}\frac {w_1(s)^2} {w_0(s)^2}\le b^2+c^2\frac {b^2+r^2c^2}
{(b^2+rc^2)^2}< b^2+1.$$ Similarly, (\ref{upl}) yields the
estimate
$$\frac {g_0^2(s)+rg_1^2(s)} {w_1^2(s)}=b^2\frac {w_0^2(s)}
{w_1^2(s)}+\frac {rc^2} {b^2+rc^2}\le b^2(b^2+c^2/r)+\frac {rc^2}
{b^2+rc^2}< b^2+1.$$ We conclude that $c_a<\sqrt{b^2+1}$. It
follows that $\gamma(\vec{G})<\sqrt{2}$.
\end{proof}

\begin{remark} The above proposition combined with a simple application of
Proposition \ref{gcont} and John's theorem, and also with Theorem
\ref{sprop}, shows that $\gamma(\vec{X})<2$ for every
two-dimensional (real) Banach couple $\vec{X}$.

\end{remark}

\subsubsection{\label{discussion}Further discussion.}
From here onwards, in view of Remark {\ref{4lat}, and since we have
seen that $c_0=c_{\pi/2}=1$, we need only consider the case where
$a\in (0,\pi/2)$ and so the numbers $b$ and $c$ are strictly
positive.

We will also suppose that $r>1$ (cf. Remark \ref{rwlog}).

Let $T=T_{g_0,g_1}$ be a member of $\mathcal{T}_{b,c}$ for which
the infimum
\begin{equation}
c_a=\inf_{T\in \mathcal{T}_{b,c}}\left\| T\right\|
_{\vec{P}\rightarrow \vec{G}}. \label{ninf}
\end{equation}
is attained. Lemma \ref{findim} guarantees that such an operator
$T$ exists.

The exact value of $c_a$ evades us at this point, but we hope that
the following remarks will provide a step on the way towards
calculating $c_a$ and therefore also $\gamma(\vec{G})$. We will show
below that the functions $g_0$, $g_1$ possess certain properties. We
will also prove the estimate $c_a<(1+\sqrt{r})/2$. This will imply,
in view of Proposition \ref{sq2prop}, that
\begin{equation}\label{ggr}
\gamma(\vec{G})<\min\left\{ \frac {1+\sqrt{r}} 2,\sqrt{2} \right\}.
\end{equation}

\begin{remark}\label{4flt2} Let $\widetilde{g}_{0}$ and $\widetilde{g}_{1}$ be the functions defined by $%
\widetilde{g}_{j}:=\frac{\left| \int_{0}^{\infty}g_{j}(s )ds \right| }{%
\int_{0}^{\infty}\left| g_{j}(s )\right| ds }|g_{j}|$ for $j=0,1$.
It is easy
to check that the operator $\widetilde{T}=T_{\widetilde{g}_{0},\widetilde{g}%
_{1}}$ is also in $\mathcal{T}_{b,c}$ and that $\|T_{\widetilde{g}_{0},\widetilde{g}%
_{1}}\|_{\vec{P}\to\vec{G}}\le\|T_{g_0,g_1}\|_{\vec{P}\to\vec{G}}$.
\end{remark}

By Remark \ref{4flt2}, we can and will assume that $g_0$ and $g_1$
are non-negative a.e. The conditions on $T$ imply that
\begin{equation}
\left\{
\begin{array}{lll}
g_{0}^{2}+g_{1}^{2} & \le c_a^{2}w_{0}^{2} & \text{and also} \\
g_{0}^{2}+rg_{1}^{2} & \le c_a^{2}w_{1}^{2} & \text{at almost every point of }%
(0,\infty )\text{.}
\end{array}
\right. .  \label{japc}
\end{equation}
We introduce two subsets $E_0$, $E_1$ of $(0,\infty)$ defined by
$$E_i=\{s\in (0,\infty):g_0(s)^2+r^i
g_1(s)^2=c_a^2w_i(s)^2\},\quad i=0,1.$$ The following simple fact
is true.

\begin{fact} \label{zjapc} The set $E_0\cup E_1$ contains almost every point of $(0,\infty)$.
\end{fact}

\begin{proof} Suppose, on the contrary, that there exists a set $%
E\subset (0,\infty )$ of positive measure, such that $%
g_{0}^{2}+g_{1}^{2}<c_a^{2}w_{0}^{2}$ and also $%
g_{0}^{2}+rg_{1}^{2}<c_a^{2}w_{1}^{2}$ at every point of $E$. Then
we can suppose, replacing $E$ if necessary by a smaller subset
also having positive measure, that, for some positive $\epsilon $
,
\begin{equation}
g_{0}^{2}+g_{1}^{2}<(1-\epsilon )c_a^{2}w_{0}^{2}\,\text{and also }%
g_{0}^{2}+rg_{1}^{2}<(1-\epsilon )c_a^{2}w_{1}^{2}\text{ at all
points of }E. \label{fluq}
\end{equation}

\smallskip For $j=0,1$ we define the function $\widetilde{g}_{j}=\sqrt{%
g_{j}^{2}+\phi }$ where
\begin{equation}
\phi =\epsilon c_a^{2}\chi _{E}\min \left\{ \frac{w_{0}^{2}}{2},\frac{w_{1}^{2}%
}{1+r}\right\} .  \label{zaml}
\end{equation}
It follows easily from (\ref{japc}), (\ref{fluq}) and (\ref{zaml})
that, for $j=0,1$, we have
\begin{equation}
\widetilde{g}_{0}^{2}+r^{j}\widetilde{g}%
_{1}^{2}=g_{0}^{2}+r^{j}g_{1}^{2}+(1+r^{j})\phi \le
c_a^{2}w_{j}^{2} \label{bllamp}
\end{equation}
at every point of $E$ and at almost every point of $(0,\infty
)\backslash E$.

Since $w_{0}$ and $w_{1}$ are both strictly positive on $(0,\infty
)$ and $E$ has positive measure, it follows that
\begin{equation}
\widetilde{b}:=\int_{0}^{\infty }\widetilde{g}_{0}(s)ds>b=\text{ }%
\int_{0}^{\infty }g_{0}(s)ds\text{ and }\widetilde{c}:=\int_{0}^{\infty }%
\widetilde{g}_{1}(s)ds>\int_{0}^{\infty }g_{1}(s)ds=c \label{broj}
\end{equation}
and so the operator $S$ defined by $S=T_{v_{0},v_{1}}$ where $v_{0}=\frac{b}{%
\widetilde{b}}\widetilde{g}_{0}$ and $v_{1}=\frac{c}{\widetilde{c}}%
\widetilde{g}_{1}$ satisfies $Sf=\alpha $. In view of
(\ref{bllamp}), (\ref{broj}) and (\ref{oinq}), its norm satisfies $\left\| S\right\| _{\vec{P}%
\rightarrow \vec{G}}\le \max \left\{ \frac{b}{\widetilde{b}},\frac{c}{%
\widetilde{c}}\right\} c_a<c_a$. This contradicts the minimal
property of $c_a$, i.e.\ (\ref{ninf}), and so proves
(\ref{zjapc}).
\end{proof}

It is convenient to
restate
Fact \ref{zjapc}
slightly differently as:
\begin{equation*}
\text{For a.e.\ }s\in (0,\infty )\text{ the point
}(g_{0}(s),g_{1}(s))\in \partial
\newEE _{s},
\end{equation*}
where the sets $\newEE _{s}$ are defined by
\begin{equation*}
\newEE _{s}\smallskip :=\left\{ (x,y):x\ge 0,y\ge 0,x^{2}+y^{2}\le
c_a^{2}w_{0}^2(s),x^{2}+ry^{2}\le c_a^{2}w_{1}^{2}(s)\right\} .
\end{equation*}

The boundary of $\newEE _{s}$ consists of a segment of the $x$ axis, a
segment of the $y$ axis, and subsets of the quarter circle $C_{s}$
of radius $c_aw_{0}(s)$
and of the quarter ellipse $\Gamma _{s}$ with semi-axes of lengths $%
c_aw_{1}(s) $ and $\frac{1}{\sqrt{r}}c_aw_{1}(s)$ in the
directions of the $x$ and $y$ axes respectively.

Since $r>1$ we see from (\ref{skm}) that
\begin{equation}
w_{0}(s)<w_{1}(s)  \label{jz}
\end{equation}
and so, on and slightly above the $x$ axis, the points of $\Gamma
_{s}$ lie strictly to the right of $C_{s}$. On the other hand,
since we shall show that
\begin{equation}
w_{0}(s)>\frac{1}{\sqrt{r}}w_{1}(s),  \label{ntc}
\end{equation}
it will follow that the points of $C_{s}$ on and near the $y$ axis
lie strictly above $\Gamma _{s}$. The sets $C_{s}$ and $\Gamma
_{s}$ intersect at a single point $(x(s),y(s))$ whose exact
coordinates will be calculated in a moment. In view of (\ref{jz})
and (\ref{ntc}) we will be able to assert that, apart from parts
of the $x$ and $y$ axes, the boundary of $\newEE _{s}$
consists of the circular arc $C_{s}^{*}$ of radius $c_aw_{0}(s)$ from $%
(c_aw_{0}(s),0)$ to $(x(s),y(s))$ and the portion $\Gamma
_{s}^{*}$ of the
quarter ellipse $\Gamma _{s}$ from $(x(s),y(s))$ to $(0,\frac{1}{\sqrt{r}}%
c_aw_{1}(s))$.

Let us now prove (\ref{ntc}). Using (\ref {fim}) and (\ref{lol})
we see that it suffices to show that $\frac {b^2+rc^2}
{b^2+r^2c^2}>\frac{1}{r}$, which is clear, since $rb^{2}>b^{2}$.

\smallskip To obtain explicit expressions for $x(s)$ and $y(s)$ we simply
solve the two equations
\begin{equation}\label{mozart}x(s)^{2}+y(s)^{2}=c_a^{2}w_{0}(s)^{2}
\ \text{ and } \
x(s)^{2}+ry(s)^{2}=c_a^{2}w_{1}(s)^{2}
\end{equation}
which gives
$y(s)^{2}=\frac{%
c_a^{2}(w_{1}(s)^{2}-w_{0}(s)^{2})}{r-1}$ and then $x(s)^{2}=\frac{%
c_a^{2}(rw_{0}(s)^{2}-w_{1}(s)^{2})}{r-1}$. From this we deduce that
\begin{equation}\label{xeq}x(s)=c_aw_j(s)\frac {b(1+rs)}
{\sqrt{b^2(1+rs)^2+r^{1+j}c^2(1+s)^2}},\quad j=0,1,\end{equation}
and
\begin{equation}\label{yeq}
y(s)=c_aw_j(s)\frac {c\sqrt{r}(1+s)}
{\sqrt{b^2(1+s)^2+r^{1+j}c^2(1+rs)^2}},\quad j=0,1.\end{equation}

\begin{remark}\label{bad}
In addition to Fact \ref{zjapc} it is now plain that, for the
optimal functions $g_0$ and $g_1$ we have
$$g_0(s)\ge x(s) \text{ and }  g_1(s)\le y(s) \text{  on
} E_0$$ and likewise
$$g_0(s)\le x(s) \text{ and }  g_1(s)\ge y(s) \text{  on
} E_1.$$ At first glance one might suspect that
$E_0=E_1=(0,\infty)$, i.e., that $g_0(s)=x(s)$ and $g_1(s)=y(s)$.
However, if this were the case, we would have that
$$\int_0^\infty
\frac {x(s)} {bc_a}ds\ge \frac 1 {c_a}>\frac 1 {\sqrt{2}} \text{ and
}\int_0^\infty \frac {y(s)} {cc_a}ds\ge \frac 1 {c_a}>\frac 1
{\sqrt{2}},$$ where we have used Proposition \ref{sq2prop}. On the
other hand, a numerical calculation making use of the explicit
formula (\ref{xeq}) with the values $r=1000$, $b=\sqrt{3}/2$ and
$c=1/2$ yields $\int_0^\infty (x(s)/c_ab)ds\approx
0.6896<1/\sqrt{2}$. Thus the functions $x$ and $y$ are not optimal
in general.
\end{remark}

\smallskip
We shall now use the operators $T=T_{x/c_a,y/c_a}$ to obtain some
new information about $\gamma(\vec{G})$. From (\ref{mozart}) and
(\ref{oinq}) it is evident that $\|T\|_{\vec{P}\to \vec{G}}=1$. In
order to prove the estimate $c_a<(1+\sqrt{r})/2$ it clearly suffices
to prove that $(Tf)_1>\frac {2b} {1+\sqrt{r}}$ and $(Tf)_2>\frac
{2c} {1+\sqrt{r}}$, i.e.,
\begin{equation}\label{bugly}
\int_0^\infty \frac {x(s)}{bc_a}ds>\frac 2{1+\sqrt{r}} \text{ and
}\int_0^\infty \frac {y(s)} {cc_a}ds>\frac 2
{1+\sqrt{r}}.\end{equation} In order to prove (\ref{bugly}), we
observe that, for $j=0,1$, the functions
$$u_j(s):=\frac {1+rs}
{\sqrt{b^2(1+rs)^2+r^{1+j}c^2(1+s)^2}}= 1/\sqrt{b^2+r^{1+j}c^2\frac
{(1+s)^2} {(1+rs)^2}}
$$
are increasing on $(0,\infty)$ and, likewise, the functions
$$
v_j(s):=\frac {\sqrt{r}(1+s)}
{\sqrt{b^2(1+s)^2+r^{1+j}c^2(1+rs)^2}}
$$
are decreasing on
$(0,\infty)$. By (\ref{xeq}) we obtain
$$\int_0^\infty \frac {x(s)} {bc_a}ds=\int_0^{1/\sqrt{r}}u_0(s)d\gamma_0(s)+\int_{1/\sqrt{r}}^\infty
u_0(s)d\gamma_0(s)>$$
$$>
u_0(0)(\gamma_0(1/\sqrt{r})-\gamma_0(0))+u_0(1/\sqrt{r})(1-\gamma_0(1/\sqrt{r}))=$$
$$=\frac 1 {\sqrt{b^2+rc^2}}\cdot \frac {\sqrt{b^2+rc^2}}
{1+\sqrt{r}}+ 1\cdot\biggl(1-\frac {\sqrt{b^2+rc^2}}
{1+\sqrt{r}}\biggr)=$$
$$=1+\frac {1-\sqrt{b^2+rc^2}} {1+\sqrt{r}}\ge 1+\frac {1-\sqrt{r}} {1+\sqrt{r}}=
\frac 2 {1+\sqrt{r}}.$$ Similarly, by using (\ref{yeq}), we get
$$\int_0^\infty \frac {y(s)} {cc_a} ds>
v_1(1/\sqrt{r})(\gamma_1(0)-\gamma_1(1/\sqrt{r}))+v_1(\infty)\gamma_1(1/\sqrt{r})=$$
$$=\frac 1 {\sqrt{b^2+rc^2}}\cdot \biggl(\sqrt {b^2+rc^2}-\frac
{\sqrt{r}} {1+\sqrt{r}}\biggr)+\frac {\sqrt{r}}{1+\sqrt{r}}\cdot
\frac 1 {\sqrt{r}}\ge $$
$$\ge 1+\frac {\sqrt{r}} {1+\sqrt{r}}\biggl(\frac 1 {\sqrt{r}}-1\biggr)=1+\frac
{1-\sqrt{r}} {1+\sqrt{r}}=\frac 2 {1+\sqrt{r}}.
$$
This establishes (\ref{bugly}) and so indeed we have
$c_a<(1+\sqrt{r})/2$ and can deduce (\ref{ggr}).

\smallskip

\section{\label{latticecouple}\smallskip The two dimensional couple
$\vec{X}=(\ell ^{2}_2,\ell ^{\infty }_2)$}

\subsection{\label{terminology}Terminology, notation and some
preliminaries.} Let $\vec{X}=(\ell^2_2,\ell^\infty_2)$. Consider the
point $\alpha =(1,a)\in X_{0}+X_{1}$ where $a>1$.

Let $E(t,\alpha ;\vec{X})$ be the error functional $$E(t,\alpha ;\vec{X}%
)=\inf \left\{ \left\| \alpha -\beta \right\| _{X_{0}}:\beta \in
X_{1},\left\| \beta \right\| _{X_{1}}\le t\right\}.$$ Then, for $t\in (0,1]$%
, the optimal choice of $\beta $ is $(t,t)$. For $t\in [1,a]$ the
optimal choice of $\beta $ is $(t,1)$, and for $t>a$ the optimal
choice is $\beta =\alpha $. Consequently
$$E(t,\alpha
;\vec{X})=\left\{
\begin{array}{lll}
\sqrt{(a-t)^{2}+(1-t)^{2}} & , & 0\le t\le 1 \\
a-t & , & 1<t\le a \\
0 & , & t>a
\end{array}
\right. .$$

Now let $w:(0,a)\rightarrow (1,\infty )$ be a non increasing
function and consider the couple of weighted $L^{1}$ spaces
$\vec{P}=(P_{0},P_{1})$ on
the measure space $(0,a)$ (equipped with Lebesgue measure) where $%
P_{0}=L_{w}^{1}$ and $P_{1}=L^{1}$. Let $f=\chi _{(0,a)}$, and let $E(t,f;%
\vec{P})=\inf \left\{ \left\| f-g\right\| _{P_{0}}:g\in
P_{0},\left\| g\right\| _{P_{1}}\le t\right\} $. Since $w\ge 1$
and $w$ is non increasing, the optimal choice for $g$ is $\chi
_{[0,\min (t,a)]}$ for all $t\in (0,\infty )$. It follows that
$E(t,f;\vec{P})=\left\| f-g\right\| _{P_{0}}=\int_{\min
(t,a)}^{a}w(\xi )d\xi $.

If $w$ is continuous, then $E(t,f;\vec{P})$ is differentiable,
with derivative equal to $-w(t)$ for all $t\in (0,a)$.

The function $E(t,\alpha ;\vec{X})$ is also differentiable on
$(0,a)$ and its derivative for $t\in (0,a)$ is given by
\begin{equation*}
\frac{d}{dt}E(t,\alpha ;\vec{X})=\left\{
\begin{array}{lll}
\frac{2t-a-1}{\sqrt{(a-t)^{2}+(1-t)^{2}}} & , & 0<t<1 \\
-1 & , & 1\le t<a
\end{array}
\right. .
\end{equation*}

\smallskip By general properties of the error functional, this derivative
must be negative and non-decreasing. Thus the function
\begin{equation}
w_{*}(t):=-\frac{d}{dt}E(t,\alpha ;\vec{X})=\left\{
\begin{array}{lll}
\frac{a+1-2t}{\sqrt{(t-a)^{2}+(t-1)^{2}}} & , & 0<t<1 \\
1 & , & 1\le t<a
\end{array}
\right.  \label{defw}
\end{equation}
is continuous and non-increasing and $w_{*}(t)\ge 1$ on $(0,a)$.
In fact, as can be shown directly, it is strictly decreasing on
$(0,1]$. If we now choose $w=w_{*}$ then it is easy to check that
$E(t,f;\vec{P})=$ $E(t,\alpha ;\vec{X})$ for all $t>0$. This is
equivalent, using well known connections between error
functionals, $K$-functionals and the Gagliardo diagram, to the
condition
\begin{equation}
K(t,f;\vec{P})=K(t,\alpha ;\vec{X})\text{ for all }t>0.
\label{kfe}
\end{equation}

For the rest of this section $w$ will always denote the particular
function defined by (\ref{defw}), for some choice of the constant
$a$. It is easy to check that, for every choice of $a>1$, we have
\begin{equation}
1\le w(t)<\sqrt{2}\text{ , and so also }\sqrt{w^{2}(t)-1}<1\text{
, for all }t\in (0,a)  \label{wt}
\end{equation}

\smallskip For each fixed $a\ge 1$, let $\mathcal{T}_{a}$ be the set of all
bounded linear operators $T:\vec{P}\rightarrow \vec{X}$, which,
for $f=\chi _{(0,a)}$ and $\alpha =(a,1)$ and $w$ as above,
satisfy $Tf=\alpha $.

\smallskip Let $T$ be an arbitrary operator in $\mathcal{T}_{a}$. Then $T$
has the form
$$Th=\left( \lambda _{0}(h),\lambda _{1}(h)\right)
\text{ for all }h\in P_{0}+P_{1},$$ where $\lambda _{0}$ and
$\lambda _{1}$ are both elements of $\left( P_{0}\right) ^{*}\cap
(P_{1})^{*}$ such that
$$\lambda _{0}(\chi _{(0,a)})=a\text{ and
}\lambda _{1}(\chi _{(0,a)})=1.$$
The norm of $T$ satisfies $\left\| T\right\| _{\vec{P}\rightarrow \vec{X}%
}\le c$ for some positive constant $c$, if and only if
$$\left\| \lambda _{j}\right\| _{\left( P_{0}\right) ^{*}}\le c\text{ for }%
j=0,1$$ and
$$\left| \lambda _{0}(h)\right| ^{2}+\left| \lambda
_{1}(h)\right| ^{2}\le c^{2}\left\| h\right\| _{P_{1}}^{2}\text{
for all }h\in P_{1}.$$

\smallskip We are interested in the quantity
\begin{equation}
c_{a}:=\inf \left\{ \left\| T\right\| _{\vec{P}\rightarrow
\vec{X}}:T\in \mathcal{T}_{a}\right\} .  \label{tinf}
\end{equation}
By (\ref{kfe}) and standard properties of the $K$-functional we
clearly have that
\begin{equation}
c_{a}\ge 1.  \label{stum}
\end{equation}

By Lemma \ref{findim} the infimum in (\ref{tinf}) is attained for
some $T\in \mathcal{T}_{a}$.

\smallskip There is of course a more concrete version of the representation
given above for operators $T\in \mathcal{T}_{a}$:

\smallskip In general, every bounded linear operator $T:\vec{P}\rightarrow
\vec{X}$ is determined by two functions $g_{0}$ and $g_{1}$ in
$L^{\infty }(0,a)$. More specifically we will use the notation
$T=T_{g_{0},g_{1}}$, where
\begin{equation}
Th=T_{g_{0},g_{1}}h=\left( \int_{0}^{a}h(\xi )g_{0}(\xi )d\xi
,\int_{0}^{a}h(\xi )g_{1}(\xi )d\xi \right) \text{ for each }h\in
P_{0}+P_{1} \label{glunk}
\end{equation}

Such an operator $T_{g_{0},g_{1}}$ is in $\mathcal{T}_{a}$ if and
only if the functions $g_{0}$ and $g_{1}$ also satisfy
\begin{equation}
\int_{0}^{a}g_{0}(\xi )d\xi =a\text{ and }\int_{0}^{a}g_{1}(\xi
)d\xi =1. \label{oonk}
\end{equation}

For any $T_{g_{0},g_{1}}:\vec{P}\rightarrow \vec{X}$, the norm estimate $%
\left\| T_{g_{0},g_{1}}\right\| _{\vec{P}\rightarrow \vec{X}}\le
c$ is equivalent to the two conditions
\begin{equation}
\left\| g_{j}\right\| _{L^{\infty }}\le c\text{ for }j=0,1
\label{litg}
\end{equation}
and
\begin{equation}
\left( \int_{0}^{a}h(\xi )g_{0}(\xi )d\xi \right) ^{2}+\left(
\int_{0}^{a}h(\xi )g_{1}(\xi )d\xi \right) ^{2}\le c^{2}\left(
\int_{0}^{a}|h(\xi )|w(\xi )d\xi \right) ^{2}\text{ for all }h\in
P_{1}. \label{jse}
\end{equation}
In fact (\ref{jse}) is equivalent to
\begin{equation}
g_{0}(\xi )^{2}+g_{1}(\xi )^{2}\le c^{2}w^{2}(\xi )\text{ for
a.e.\ }\xi \in (0,a).  \label{zjse}
\end{equation}
The proof that (\ref{jse}) implies (\ref{zjse}) follows readily
from the Lebesgue differentiation theorem. The reverse implication
follows easily from a suitable version of Minkowski's inequality
or Schwartz' inequality.

\subsection{\label{simple}A simple estimate from below for $\gamma
(\ell _{2}^{2},\ell _{2}^{\infty })$.}

\smallskip We can now easily show that $\vec{X}=(\ell _{2}^{2},\ell
_{2}^{\infty })$ is an example, perhaps the simplest known example
so far, of a Banach couple whose $K$-divisibility constant
satisfies
\begin{equation}
\gamma (\vec{X})>1.  \label{bto}
\end{equation}
It will be convenient to use the terminology \textit{\textbf{not exactly}}%
\textbf{\ }$K$\textbf{-}\textit{\textbf{divisible} (n.e.K-d.)} for
any Banach couple satisfying (\ref{bto}).

\begin{remark}
This example is of interest for a number of reasons:

$\bullet $ It is apparently the first known example of a couple of
rearrangement invariant spaces which is n.e.K-d.

$\bullet $ It also shows that there is no ``tight'' connection
between the exact $K$-divisibility property and the exact
Calder\'{o}n property. Neither of the couples $(\ell _{2}^{2},\ell
_{2}^{\infty })$ and $(\ell _{8}^{2},\ell _{8}^{\infty })$ are
exactly $K$-divisible. However, as shown in the appendix (Section
\ref{appendix}), $(\ell _{2}^{2},\ell _{2}^{\infty })$ is an exact
Calder\'{o}n couple, but $(\ell _{8}^{2},\ell _{8}^{\infty }) $ is
not. \textit{(Meanwhile we also know (see Section \ref{hc}) that
the
Hilbert couple }$\vec{Y}=\left( \ell _{2}^{2},\ell _{1}^{2}\right) $\textit{%
\ which is an exact Calder\'{o}n couple (see \cite{ameur}) is not
exactly }$K $\textit{-divisible.)} There are also exactly
$K$-divisible couples which are not exact Calder\'{o}n couples, or
not even Calder\'{o}n couples, an example is provided by the
couple $(L^{1}\oplus L^{\infty },L^{\infty }\oplus L^{1})$.

$\bullet $ We can also now see that there is not a ``tight''
connection between exact $K$-divisibility and the property of
exact monotonicity,
introduced and studied in \cite{ckeich}. This follows from the fact that $%
(\ell _{2}^{2},\ell _{2}^{\infty })$ and $(\ell _{8}^{2},\ell
_{8}^{\infty }) $ are both exactly monotone (See \cite{ckeich}
Theorem 2.1, p.\ 32). A connection between the $K$-divisibility
and monotonicity constants of a couple was established in
\cite{ckeich} (See formula (52) on page 55 of \cite {ckeich}.)
This result was strengthened in \cite{cw1}.
\end{remark}

Using well known results concerning $K$-divisibility (Theorem
\ref{ec}) it is easy to see that
$$\gamma (\vec{X})=\sup_{a\ge 1}c_{a}$$
 where $c_{a}$ is defined by
(\ref{tinf}). We shall show that $c_{a}>1$ for every $a>1$.

Suppose, on the contrary, that $c_{a}=1$ for some $a>1$. (Recall \ref{stum}%
).) Let $T$ be the operator in $\mathcal{T}_{a}$ whose existence
we established above, which satisfies $\left\| T\right\|
_{\vec{P}\rightarrow
\vec{X}}=c_{a}=1$. Then there exist functions $g_{0}$ and $g_{1}$ in $%
L^{\infty }(0,a)$ satisfying (\ref{oonk}) and also satisfying the estimates (%
\ref{litg}) and (\ref{zjse}) for $c=1$. In particular, since $%
\int_{0}^{a}g_{0}(\xi )d\xi =a$ and $\left| g_{0}(\xi )\right| \le
1$ for a.e.\ $\xi \in (0,a)$, we must have $g_{0}(\xi )=1$ a.e. It
follows that
\begin{equation}
\begin{split}
1& =\int_{0}^{a}g_{0}(\xi )d\xi \le \int_{0}^{a}\sqrt{w^{2}(\xi
)-1}d\xi =\int_{0}^{1}\sqrt{w^{2}(\xi )-1}d\xi =\\
&=\int_{0}^{1}\sqrt{\frac{(2\xi -a-1)^{2}}{(a-\xi )^{2}+(1-\xi
)^{2}}-1}d\xi .\\  \label{gac}
\end{split}
\end{equation}

\smallskip The expression under the square root in the last integral can be
rewritten as

\begin{eqnarray}
\frac{(a+1-2\xi )^{2}-(a-\xi )^{2}-(1-\xi )^{2}}{(a-\xi
)^{2}+(1-\xi )^{2}}
&=&\frac{\left( (a-\xi )+(1-\xi )\right) ^{2}-(a-\xi )^{2}-(1-\xi )^{2}}{%
(a-\xi )^{2}+(1-\xi )^{2}}  \notag \\
&=&\frac{2(a-\xi )(1-\xi )}{(a-\xi )^{2}+(1-\xi )^{2}}.
\label{agac}
\end{eqnarray}
This equals $1$ for all $\xi $ if $a=1$. But, for all $a>1$, we have $\frac{%
2(a-\xi )(1-\xi )}{(a-\xi )^{2}+(1-\xi )^{2}}<1$ for all $\xi $.
This shows that (\ref{gac}) cannot hold, and so provides the
contradiction which proves that $c_{a}>1$ and also establishes
(\ref{bto}).

\begin{remark}
\label{aeo}It is easy to show that $c_{a}=1$ when $a=1$. In this
case the function $w$ assumes the constant value $\sqrt{2}$ on
$(0,a)=(0,1)$ and the operator $T=T_{g_{0},g_{1}}$, which is
obtained by simply choosing $g_{0}$ and $g_{1}$ to be both
identically $1$, is in $\mathcal{T}_{a}$ and satisfies $\left\|
T\right\| _{\vec{P}\rightarrow \vec{X}}=1$.
\end{remark}

\subsection{\label{elaborate}A more elaborate calculation.}

\smallskip Throughout this section $a$ will denote a fixed number satisfying
$a>1$, and $g_{0}$ and $g_{1}$ will denote two particular functions in $%
L^{\infty }(0,a)$ which satisfy (\ref{glunk}) and (\ref{oonk}) for
an
operator $T_{g_{0},g_{1}}\in \mathcal{T}_{a}$ which attains the infimum $%
c_{a}$ in (\ref{tinf}). Therefore $g_{0}$ and $g_{1}$ satisfy
(\ref{litg}) and (\ref{zjse}) with $c=c_{a}$. Our goal here will
be to show that $g_{0}$ and $g_{1}$ necessarily have certain
properties. Our calculations in this section will also lead to the
estimate $\gamma(\vec{X})\le \frac{4+3\sqrt{2}}{4+2\sqrt{2}}$.

By familiar arguments (cf. Remark \ref{4flt2}) we can and will
assume that $g_{0}$ and $g_{1} $ are both non negative.

\smallskip We will use the following very simple claim several times in
subsequent steps of our argument:

\begin{claim}
\label{uscon}Suppose that $\widetilde{g}_{0}$ and
$\widetilde{g}_{1}$ are two non-negative functions in $L^{\infty
}(0,a)$ which satisfy
\begin{equation}
\int_{0}^{a}\widetilde{g}_{0}(\xi )d\xi >a\text{ and }\int_{0}^{a}\widetilde{%
g}_{1}(\xi )d\xi >1.  \label{isms}
\end{equation}
Then
\begin{equation*}
\left\| T_{\widetilde{g}_{0},\widetilde{g}_{1}}\right\|
_{\vec{P}\rightarrow \vec{X}}>c_{a}.
\end{equation*}
\end{claim}

\textit{Proof.} Suppose, on the contrary that
\begin{equation}
\left\| T_{\widetilde{g}_{0},\widetilde{g}_{1}}\right\|
_{\vec{P}\rightarrow \vec{X}}\le c_{a}.  \label{seps}
\end{equation}
Then the operator $S$ defined by
\begin{equation*}
Sh=\left( \frac{a}{\int_{0}^{a}\widetilde{g}_{0}(\xi )d\xi }\int_{0}^{a}%
\widetilde{g}_{0}(\xi )h(\xi )d\xi ,\frac{1}{\int_{0}^{a}\widetilde{g}%
_{1}(\xi )d\xi }\int_{0}^{a}\widetilde{g}_{1}(\xi )h(\xi )d\xi
\right)
\end{equation*}
has norm $\left\| S\right\| _{\vec{P}\rightarrow \vec{X}}$
strictly smaller than $c_{a}$ . But $S\in \mathcal{T}_{a}$ and so
we have a contradiction, which proves the claim.\smallskip \qed

It will be convenient to define the planar set
$$E_{\xi
}=\{(x,y)\in \Bbb{R}^{2}:0\le x\le c_{a},0\le y\le
c_{a},x^{2}+y^{2}\le c_{a}^{2}w^{2}(\xi )\}$$
for each $\xi \in
(0,a)$. Then, reformulating our remarks above, for any non
negative measurable functions $u_{0}$ and $u_{1}$ on $(0,a)$,
$\left\| T_{u_{0},u_{1}}\right\| _{\vec{P}\rightarrow \vec{X}}\le
c_{a}$ if and only if $\left( u_{0}(\xi ),u_{1}(\xi )\right) \in
E_{\xi }$ for a.e.\ $\xi \in (0,a)$. In particular, the two
particular norm minimizing functions $g_{0}$ and $g_{1}$ which we
are studying, satisfy this condition.

We note that the boundary of $E_{\xi }$ consists of two horizontal
and two vertical line segments and a circular arc of radius
$c_{a}w(\xi )$ which we will denote by $\Gamma _{\xi }$. We let
$V_{\xi }$ denote the vertical segment of the right side of the
boundary of $E_{\xi ,}$ i.e.\
$$V_{\xi }=\left\{ (c_{a},y):0\le
y\le c_{a}\sqrt{w^{2}(\xi )-1}\right\}$$
The uppermost point of
$V_{\xi }$, which is also the lowest point of $\Gamma _{\xi }$, is
\begin{eqnarray}
\left( c_{a},c_{a}\sqrt{w^{2}(\xi )-1}\right) &=&\left( c_{a}w(\xi
)\cos
\psi (\xi ),c_{a}w(\xi )\sin \psi (\xi )\right)  \notag \\
\text{where }\psi (\xi ) &=&\arctan \sqrt{w^{2}(\xi )-1}=\arccos \frac{1}{%
w(\xi )}  \label{defpsi}
\end{eqnarray}

\ Let $\mathcal{U}_{a}$ be the family of all couples
$(u_{0},u_{1})$ of non negative functions in $L^{\infty }(0,a)$
which satisfy

(i) $(u_{0}(\xi ),u_{1}(\xi ))\ne (0,0)$ for a.e.\ $\xi \in
(0,a)$, and

(ii) $\left\| T_{u_{0},u_{1}}\right\| _{\vec{P}\rightarrow
\vec{X}}\le c_{a}$ or, equivalently $\left( u_{0}(\xi ),u_{1}(\xi
)\right) \in E_{\xi }$ for a.e.\ $\xi \in (0,a)$.

We claim that the special functions $g_{0}$ and $g_{1}$ satisfy
\begin{equation}
(g_{0},g_{1})\in \mathcal{U}_{a}.  \label{rbz}
\end{equation}
They of course satisfy part (ii) of the definition. To show that
they also satisfy part (i), let
\begin{equation*}
N=\left\{ \xi \in (0,a):\left( g_{0}(\xi ),g_{1}(\xi )\right)
=(0,0)\right\}
\end{equation*}
and let $\widetilde{g}_{j}=g_{j}\chi _{(0,a)\backslash N}+\frac{c_{a}}{\sqrt{%
2}}w\chi _{N}$\ for $j=0,1$. In view of (\ref{wt}) it is clear that $(%
\widetilde{g}_{0}(\xi ),\widetilde{g}_{1}(\xi ))\in E_{\xi }$\ for a.e.\ $%
\xi \in (0,a)$, which is equivalent to (\ref{seps}). But, if $N$\
has positive measure, then (\ref{isms}) also holds, which, by
Claim \ref{uscon}, is impossible.

It is convenient to represent each $(u_{0},u_{1})\in
\mathcal{U}_{a}$ in the
``polar'' form $(u_{0},u_{1})=(\rho \cos \theta ,\rho \sin \theta )$ where $%
\rho :(0,a)\rightarrow (0,\sqrt{2})$ and $\theta :(0,a)\rightarrow [0,\frac{%
\pi }{2}]$ are the measurable functions defined by $\rho (\xi )=\sqrt{%
u_{0}^{2}(\xi )+u_{1}^{2}(\xi )}$ and $\theta (\xi )=\arcsin
\frac{u_{1}(\xi
)}{\rho (\xi )}$ for all $\xi \in (0,a)$. Accordingly, we let $\mathcal{P}%
_{a}$ be the family of all couples $(\rho ,\theta )$ of functions
$\rho
:(0,a)\rightarrow [0,\sqrt{2})$ and $\theta :(0,a)\rightarrow [0,\frac{\pi }{%
2}]$ such that $(\rho \cos \theta ,\rho \sin \theta )\in
\mathcal{U}_{a}$.

\smallskip

\begin{claim}
\label{sym}If $(\rho ,\theta )\in \mathcal{P}_{a}$ and $\phi
:(0,a)\rightarrow [0,\frac{\pi }{2}]$ is a measurable function
satisfying
\begin{equation*}
\theta (\xi )\le \phi (\theta )\le \frac{\pi }{4}\text{ or }\theta
(\xi )\ge \phi (\theta )\ge \frac{\pi }{4}\text{ }
\end{equation*}

for a.e.\ $\xi \in (0,a)$, then $(\rho ,\phi )\in
\mathcal{P}_{a}$.
\end{claim}

This is obvious, in view of the form of the sets $E_{\xi }$. \qed

We have now the following simple ``variational principle'':

\begin{lemma}
\label{var}Suppose that the functions $\rho $ and $\theta $
satisfy
\begin{equation}
(\rho ,\theta )\in \mathcal{P}_{a}\text{ and }g_{0}=\rho \cos
\theta \text{ and }g_{1}=\rho \sin \theta .  \label{gpol}
\end{equation}
Suppose that $A$ and $B$ are each measurable subsets of $(0,a)$
with positive measure. Suppose that $p$ , $q$ are real constants
such that, for some $\delta >0$ and each constant $t\in [0,\delta
]$, the function $\phi _{t}=\theta +tp\chi _{A}+tq\chi _{B}$
satisfies
\begin{equation}
(\rho ,\phi _{t})\in \mathcal{P}_{a}.  \label{vsep}
\end{equation}
Then at least one of the following two inequalities
\begin{equation}
p\int_{A}\rho (\xi )\sin \theta (\xi )d\xi +q\int_{B}\rho (\xi
)\sin \theta (\xi )d\xi \ge 0  \label{jkk}
\end{equation}
and
\begin{equation}
p\int_{A}\rho (\xi )\cos \theta (\xi )d\xi +q\int_{B}\rho (\xi
)\cos \theta (\xi )d\xi \le 0  \label{zjkk}
\end{equation}
must hold.
\end{lemma}

\textit{Proof.} Define $G_{0}(t)=\int_{0}^{a}\rho (\xi )\cos \phi
_{t}(\xi )d\xi $ and $G_{1}(t)=\int_{0}^{a}\rho (\xi )\sin \phi
_{t}(\xi )d\xi $ for all $t\in \Bbb{R}$. Standard arguments (e.g.\
via dominated convergence) show that $G_{0}$ and $G_{1}$ are
differentiable for all $t\in \Bbb{R}$ and their derivatives are
continuous functions of $t$ given by
\begin{equation*}
G_{0}^{\prime }(t)=-p\int_{A}\rho (\xi )\sin \phi _{t}(\xi )d\xi
-q\int_{B}\rho (\xi )\sin \phi _{t}(\xi )d\xi
\end{equation*}
and
\begin{equation*}
G_{1}^{\prime }(t)=p\int_{A}\rho (\xi )\cos \phi _{t}(\xi )d\xi
+q\int_{B}\rho (\xi )\cos \phi _{t}(\xi )d\xi .
\end{equation*}
Suppose that neither of (\ref{jkk}) and (\ref{zjkk}) hold. Then $%
G_{0}^{\prime }(0)$ and $G_{1}^{\prime }(0)$ are both strictly
positive. Thus $G_{0}$ and $G_{1}$ are both increasing functions
in some neighbourhood
of $0$. So, for some $\delta ^{\prime }\in (0,\delta ]$, we have $%
G_{0}(\delta ^{\prime })>G_{0}(0)$ and $G_{1}(\delta ^{\prime
})>G_{1}(0)$, or, in other words, the functions
$\widetilde{g}_{0}:=\rho \cos \phi _{\delta ^{\prime }}$ and
$\widetilde{g}_{1}:=\rho \sin \phi _{\delta ^{\prime }}$ satisfy
(\ref{isms}). But, in view of (\ref{vsep}), these same two
functions also satisfy (\ref{seps}). By Claim \ref{uscon} this is
impossible, so at least one of (\ref{jkk}) and (\ref{zjkk}) must
hold. \qed

\smallskip As our first application of Lemma \ref{var} we will prove that
\begin{equation}
g_{0}(\xi )\ge g_{1}(\xi )\text{ for a.e.\ }\xi \in (0,a).
\label{bst}
\end{equation}
If the functions $\rho $ and $\theta $ satisfy (\ref{gpol}) then
(\ref{bst}) is equivalent to
\begin{equation}
\theta (\xi )\le \frac{\pi }{4}\text{ for a.e.\ }\in (0,a).
\label{qbst}
\end{equation}
So, if (\ref{bst}) is false, then the set $\{\xi \in
(0,a):g_{0}(\xi )<g_{1}(\xi )\}=\{\xi \in (0,a):\theta (\xi
)>\frac{\pi }{4}\}$ has positive
measure and, furthermore, for some positive number $\eta _{0}$, the set $%
A:=\{\xi \in (0,a):\theta (\xi )>\eta _{0}+\frac{\pi }{4}\}$ also
has positive measure. Since $\int_{0}^{a}g_{0}(\xi )d\xi
>\int_{0}^{a}g_{1}(\xi )d\xi $ the set $\{\xi \in (0,a):g_{0}(\xi
)>g_{1}(\xi )\}=\{\xi \in (0,a):\theta (\xi )<\frac{\pi }{4}\}$
must also have positive measure, and so, for some positive number
$\eta _{1}$, the set $B=\{\xi \in (0,a):\theta (\xi )<\frac{\pi
}{4}-\eta _{1}\}$ also has positive measure. Let $p$ be an
arbitrary negative number and let $q=1$. Let us also choose
$\delta =\min \{\eta _{0}/|p|,\eta _{1}\}$. Then, using Claim
\ref{sym}, we see that all the hypotheses of Lemma \ref{var} hold.
Consequently, Lemma \ref{var} implies that
\begin{displaymath}
\begin{split}
& p\int_{A}\rho (\xi )\sin \theta (\xi )d\xi +\int_{B}\rho (\xi
)\sin \theta (\xi )d\xi \ge 0\text{ or }\\
&p\int_{A}\rho (\xi )\cos \theta (\xi )d\xi +\int_{B}\rho (\xi
)\cos \theta (\xi )d\xi \le 0\\
\end{split}
\end{displaymath}

But now w\smallskip e shall show that we have a contradiction by
finding a negative number $p$ which satisfies
\begin{equation}
\begin{cases}
p\int_{A}\rho (\xi )\sin \theta (\xi )d\xi +\int_{B}\rho (\xi
)\sin \theta (\xi )d\xi <0\text{ and }\cr p\int_{A}\rho (\xi )\cos
\theta (\xi )d\xi +\int_{B}\rho (\xi )\cos \theta (\xi )d\xi
>0.\cr \end{cases} \label{numi}
\end{equation}
In view of (\ref{rbz}), $\rho (\xi )>0$ for a.e.\ $\xi $ and
$\int_{A}\rho (\xi )\sin \theta (\xi )d\xi >\int_{A}\rho (\xi )\sin
\frac{\pi }{4}d\xi
>0$.
 We also have $\int_{A}\rho (\xi )\sin \theta (\xi )d\xi
>\int_{A}\rho (\xi )\cos \theta (\xi )d\xi \ge 0$. Similarly
$\int_{B}\rho (\xi )\cos \theta (\xi )d\xi >\int_{B}\rho (\xi
)\cos \frac{\pi }{4}d\xi >0$ and $\int_{B}\rho (\xi )\cos \theta
(\xi )d\xi >\int_{B}\rho (\xi )\sin \theta (\xi )d\xi \ge
0 $.

If $\int_{A}\rho (\xi )\cos \theta (\xi )d\xi =0$ then every number $p<-%
\frac{\int_{B}\rho (\xi )\sin \theta (\xi )d\xi }{\int_{A}\rho
(\xi )\sin \theta (\xi )d\xi }$ satisfies (\ref{numi}). Otherwise,
if $\int_{A}\rho (\xi )\cos \theta (\xi )d\xi \ne 0$, then
condition (\ref{numi}) is equivalent to
\begin{equation*}
p+\frac{\int_{B}\rho (\xi )\sin \theta (\xi )d\xi }{\int_{A}\rho
(\xi )\sin \theta (\xi )d\xi }<0\text{ and }p+\frac{\int_{B}\rho
(\xi )\cos \theta (\xi )d\xi }{\int_{A}\rho (\xi )\cos \theta (\xi
)d\xi }>0
\end{equation*}
and so also to
\begin{equation*}
\frac{\int_{B}\rho (\xi )\sin \theta (\xi )d\xi }{\int_{A}\rho
(\xi )\sin
\theta (\xi )d\xi }<-p<\frac{\int_{B}\rho (\xi )\cos \theta (\xi )d\xi }{%
\int_{A}\rho (\xi )\cos \theta (\xi )d\xi }\ .
\end{equation*}
So it is clear that we can find $p$ with the required properties,
if and only if
\begin{equation}
\frac{\int_{B}\rho (\xi )\sin \theta (\xi )d\xi }{\int_{B}\rho
(\xi )\cos
\theta (\xi )d\xi }<\frac{\int_{A}\rho (\xi )\sin \theta (\xi )d\xi }{%
\int_{A}\rho (\xi )\cos \theta (\xi )d\xi }.  \label{oonq}
\end{equation}
Since $\sin \theta (\xi )<\cos \theta (\xi )$ for all $\xi \in B$,
and $\sin \theta (\xi )>\cos \theta (\xi )$ for all $\xi \in A$,
the left term of (\ref {oonq}) is strictly less than $1$ and the
right term of (\ref{oonq}) is strictly greater than $1$. This
proves (\ref{oonq}) and so provides the contradiction which
establishes (\ref{bst}).

\begin{claim}
\smallskip \label{thie}For almost every $\xi \in (0,a)$, if $g_{0}(\xi
)=c_{a}$ then $g_{1}(\xi )=c_{a}\sqrt{w^{2}(\xi )-1}$ and, consequently, $%
\left( g_{0}(\xi ),g_{1}(\xi )\right) $ is the upper endpoint
$\left(
c_{a}w(\xi )\cos \psi (\xi ),c_{a}w(\xi )\sin \psi (\xi )\right) $ of $%
V_{\xi }$ as defined in (\ref{defpsi}).
\end{claim}

\textit{Proof.} This amounts to showing that the set
$$V:=\{\xi
\in (0,a):g_{0}(\xi )=c_{a},\;g_{0}^{2}(\xi )+g_{1}^{2}(\xi
)<c_{a}^{2}w^{2}(\xi )\}$$
has measure $0$. If this is not true,
then the function $u_{1}:=g_{1}\chi _{(0,1)\backslash
V}+c_{a}\sqrt{w^{2}-1}\chi _{V}$ satisfies
\begin{equation}
\int_{0}^{a}u_{1}(\xi )d\xi >\int_{0}^{a}g_{1}(\xi )d\xi =1.
\label{kubc}
\end{equation}
Furthermore (in view of (\ref{wt})) it is clear that $(g_{0}(\xi
),u_{1}(\xi
))\in E_{\xi }$ for a.e.\ $\xi \in (0,a)$. Since $c_{a}>1$ and $%
\int_{0}^{a}g_{0}(\xi )d\xi =a$ the set $V_1=\left\{ \xi \in
(0,a):g_{0}(\xi )<c_{a}\right\} $ must also have positive measure.
Let $V_{*}$ be some subset of $V_1$ which also has positive
measure and define
\begin{equation*}
\widetilde{g}_{0}=g_{0}\chi _{(0,a)\backslash V_{*}}+c_{a}\chi
_{V_{*}}\text{ and }\widetilde{g}_{1}=u_{1}\chi _{(0,a)\backslash
V_{*}}.
\end{equation*}
Then $(\widetilde{g}_{0}(\xi ),\widetilde{g}_{1}(\xi ))\in E_{\xi
}$ for a.e.\ $\xi \in (0,a)$ and
$\int_{1}^{a}\widetilde{g}_{0}(\xi )d\xi
>\int_{1}^{a}g_{0}(\xi )d\xi =a$. If we choose the measure of $V_{*}$ to be
sufficiently small then we will also have, using (\ref{kubc}), that $%
\int_{1}^{a}\widetilde{g}_{1}(\xi )d\xi >\int_{1}^{a}g_{1}(\xi
)d\xi =1$. Once again we can apply Claim \ref{uscon} to obtain a
contradiction. This proves that the set $V$ has measure $0$. \qed

\smallskip

Our next step is to show that
\begin{equation}
\text{The set }Q=\left\{ \xi \in (0,a):g_{1}(\xi )=0,\;g_{0}(\xi
)<c_{a}\right\} \text{ has measure }0.  \label{pwlm}
\end{equation}
If this is false, then we consider the functions $\widetilde{g}_{0}=\sqrt{%
\frac{1}{2}\left( g_{0}^{2}+c_{a}^{2}\right) }\chi _{Q}+g_{0}\chi
_{(0,a)\backslash Q}$ and $\widetilde{g}_{1}=\min \left\{ \sqrt{%
c_{a}^{2}w^{2}-\widetilde{g}_{0}^{2}},c_{a}\right\} \chi
_{Q}+g_{1}\chi
_{(0,a)\backslash Q}$. It is clear that on the set $Q$ we have $g_{0}<%
\widetilde{g}_{0}<c_{a}\le c_{a}w$ and consequently also $\widetilde{g}%
_{1}>0=g_{1}$. Consequently $\widetilde{g}_{0}$ and
$\widetilde{g}_{1}$
satisfy (\ref{isms}). It is also clear that $(\widetilde{g}_{0}(\xi ),%
\widetilde{g}_{1}(\xi ))\in E_{\xi }$ for a.e.\ $\xi \in (0,a)$.
We can thus
use Claim \ref{uscon} to obtain a contradiction and complete the proof of (%
\ref{pwlm}).

\begin{claim}
\label{lbytbb}Suppose that, as in Lemma \ref{var}, the functions
$\rho $ and $\theta $ satisfy (\ref{gpol}). Then
\begin{equation}
\rho (\xi )=c_{a}\min \left\{ \frac{1}{\cos \theta (\xi )},w(\xi
)\right\} \text{ for a.e.\ }\xi \in (0,a).  \label{bigrho}
\end{equation}
\end{claim}

\textit{Proof.} Let us use the notation $\widetilde{\rho }(\xi
)=c_{a}\min \left\{ \frac{1}{\cos \theta (\xi )},w(\xi )\right\}
$. In view of (\ref {qbst}), it is clear that
\begin{equation}
(\widetilde{\rho },\theta )\in \mathcal{P}_{a}  \label{sbiwhwb}
\end{equation}
and that, furthermore, $\rho (\xi )\le \widetilde{\rho }(\xi )$ for a.e.\ $%
\xi \in (0,a)$. Suppose, contrarily to what we claim, that the set
$R=\{\xi \in (0,a):\rho (\xi )<\widetilde{\rho }(\xi )\}$ has
positive measure. Let us write $R=R_{0}\cup R_{1}$ where
$R_{0}=R\cap \{\xi \in (0,a):\theta (\xi )=0\}$ and
$R_{1}=R\backslash R_{0}$. We observe that $R_{0}$ is exactly the
set $Q$ of (\ref{pwlm}) which has measure $0$. Consequently
$R_{1}$ has
positive measure. This implies that the functions $\widetilde{g}_{0}=%
\widetilde{\rho }\cos \theta $ and
$\widetilde{g}_{1}=\widetilde{g}_{1}\sin \theta $ satisfy
$\int_{0}^{a}\widetilde{g}_{j}(\xi )d\xi
>\int_{0}^{a}g_{j}(\xi )d\xi $ for $j=0,1$. In view of (\ref{sbiwhwb}) and
Claim \ref{uscon} this is impossible. \qed

\smallskip

We can now show that the functions $\rho $ and $\theta $ which
satisfy (\ref {gpol}) also satisfy
\begin{equation}
\arccos \frac{1}{w(\xi )}\le \theta (\xi )\le \frac{\pi }{4}\text{
for a.e.\ }\xi \in (0,a).  \label{psth}
\end{equation}
In view of (\ref{qbst}), we can do this by showing that the set
$$W=\left\{
\xi \in (0,a):\arccos \frac{1}{w(\xi )}>\theta (\xi )\right\}$$ has measure $%
0$. Let us first observe that, by Claim \ref{thie}, almost every
$\xi \in (0,a)$ satisfying $g_{0}(\xi )=c_{a}$ also satisfies
$\theta (\xi )=\psi (\xi )=\arccos \frac{1}{w(\xi )}$ and so is
not in $W$. On the other hand, every $\xi \in W$ satisfies
$\frac{1}{w(\xi )}<\cos \theta (\xi )$. Consequently, by
(\ref{bigrho}), $\rho (\xi )=c_{a}/\cos w(\xi )$ or, equivalently,
$g_{0}(\xi )=c_{a}$ for a.e.\ $\xi \in W$. So indeed $W$ has
measure $0$ and we have proved (\ref{psth}).

\begin{theorem}
\smallskip \label{hfmg}Suppose that $\rho $ and $\theta $ are the functions
which satisfy (\ref{gpol}). Then $\theta (\xi )$ assumes a
constant value a.e.\ on the set
\begin{equation}
U=\left\{ \xi \in (0,a):\arccos \frac{1}{w(\xi )}<\theta (\xi
)\right\} . \label{defu}
\end{equation}
\end{theorem}

\textit{Proof.} Suppose that the theorem is false. Then there
exist two
subsets $A$ and $B$ of $U$, each having positive measure, and numbers $%
\theta _{0}$ and $\theta _{1}$ such that $0\le \theta _{0}<\theta
_{1}\le \pi /4$ and
\begin{equation*}
\theta (\xi )\le \theta _{0}\text{ for all }\xi \in A\text{ and
}\theta _{1}\le\theta (\xi )\text{ for all }\xi \in B.
\end{equation*}
We can assume further that each $\xi \in B$ also satisfies $\arccos \frac{1}{%
w(\xi )}<\theta (\xi )-\delta _{0}$ for some positive number
$\delta _{0}$, since, if not $B$ can be replaced by some subset of
positive measure which does have this property. Let $p=1$ and let
$q$ be an arbitrary negative
number. Then, if $\delta =\min \left\{ \frac{\pi }{4}-\theta _{0},\frac{%
\delta _{0}}{|q|}\right\} $, all the hypotheses of Lemma \ref{var}
are satisfied.

To complete the proof we will show that, for some choice of $q<0$,
both the inequalities
\begin{equation}
\int_{A}\rho (\xi )\sin \theta (\xi )d\xi +q\int_{B}\rho (\xi
)\sin \theta (\xi )d\xi <0  \label{tss}
\end{equation}
and
\begin{equation}
\int_{A}\rho (\xi )\cos \theta (\xi )d\xi +q\int_{B}\rho (\xi
)\cos \theta (\xi )d\xi >0  \label{nnu}
\end{equation}
hold and thus we have a contradiction to the conclusion which
would follow from Lemma \ref{var}.

We recall (cf.\ (\ref{rbz})) that $\rho (\xi )>0$ for a.e.\ $\xi
\in (0,a)$. So
\begin{equation*}
\int_{B}\rho (\xi )\sin \theta (\xi )d\xi \ge \int_{B}\rho (\xi
)\sin \theta _{1}d\xi = \sin \theta _{1}\int_{B}\rho (\xi )d\xi >0
\end{equation*}
and
\begin{equation*}
\int_{B}\rho (\xi )\cos \theta (\xi )d\xi \ge \int_{B}\rho (\xi )\cos \frac{%
\pi }{4}d\xi =\frac{1}{\sqrt{2}}\int_{B}\rho (\xi )d\xi >0.
\end{equation*}
Since $\tan \theta _{0}<\tan \theta _{1}$ we have
\begin{equation*}
\frac{\sin \theta _{0}}{\sin \theta _{1}}<\frac{\cos \theta
_{0}}{\cos \theta _{1}}
\end{equation*}
and consequently the numbers
\begin{equation*}
r_{0}:=\frac{\int_{A}\rho (\xi )\sin \theta (\xi )d\xi
}{\int_{B}\rho (\xi )\sin \theta (\xi )d\xi }\text{ and
}r_{1}:=\frac{\int_{A}\rho (\xi )\cos \theta (\xi )d\xi
}{\int_{B}\rho (\xi )\cos \theta (\xi )d\xi }
\end{equation*}
satisfy

\begin{equation*}
r_{0}\le \frac{\int_{A}\rho (\xi )\sin \theta _{0}d\xi
}{\int_{B}\rho (\xi
)\sin \theta _{1}d\xi }<\frac{\int_{A}\rho (\xi )\cos \theta _{0}d\xi }{%
\int_{B}\rho (\xi )\cos \theta _{1}d\xi }\le r_{1}.
\end{equation*}
Clearly every number $q$ satisfying $r_{0}<-q<r_{1}$ is negative
and also satisfies (\ref{tss}) and (\ref{nnu}). This completes the
proof of the theorem. \qed

\smallskip

Let $\theta _{a}$ be the constant value assumed a.e.\ by $\theta
(\xi )$ on the set $U$ defined by (\ref{defu}). Then, perhaps
after altering $\rho $ and $\theta $ on sets of measure $0$, we
obtain that $U=\left\{ \xi \in
(0,a):\arccos \frac{1}{w(\xi )}<\theta _{a}\right\} $. In view of (\ref{psth}%
), $\arccos \frac{1}{w(\xi )}=\theta (\xi )$ for a.e.\ $\xi \in
(0,a)\backslash U$.

If $\theta _{a}=0$, then $U$ is empty and so $w(\xi )\cos \theta
(\xi )=1$ for a.e.\ $\xi \in (0,a)$. Consequently (cf.\
(\ref{bigrho})) $\rho (\xi )=c_{a}/\cos \theta (\xi )$ for a.e.\
$\xi \in (0,a)$ and so
\begin{equation*}
\int_{0}^{a}g_{0}(\xi )d\xi =\int_{0}^{a}\rho (\xi )\cos \theta
(\xi )d\xi =c_{a}a.
\end{equation*}
But, since $c_{a}>1$, this contradicts (\ref{oonk}). We deduce
that $\theta _{a}>0$.

At the other extreme, if $\theta _{a}\ge \arccos \frac{1}{w(0)}$
then, since $w$ is strictly decreasing on $[0,1]$, we obtain that
$U=(0,a)$ and it follows from (\ref{bigrho}) that $\rho (\xi
)=c_{a}w(\xi )$ for a.e.\ $\xi \in (0,a)$. We also have
\begin{equation*}
a=\frac{\int_{0}^{a}g_{0}(\xi )d\xi }{\int_{0}^{a}g_{1}(\xi )d\xi }=\frac{%
\int_{0}^{a}\rho (\xi )\cos \theta _{a}d\xi }{\int_{0}^{a}\rho
(\xi )\sin \theta _{a}d\xi }=\tan \theta _{a},
\end{equation*}
which implies that $\sin \theta _{a}=a/\sqrt{a^{2}+1}$.
Consequently,
\begin{equation}
\int_{0}^{a}g_{1}(\xi )d\xi =\int_{0}^{a}\rho (\xi )\sin \theta _{a}d\xi =%
\frac{a}{\sqrt{a^{2}+1}}\int_{0}^{a}c_{a}w(\xi )d\xi .
\label{rnd}
\end{equation}
In view of (\ref{defw}),
\begin{equation*}
\int_{0}^{a}w(\xi )d\xi =-\int_{0}^{a}\frac{d}{d\xi }E(\xi ,\alpha ;\vec{X}%
)d\xi =E(0,\alpha ;\vec{X})-E(a,\alpha ;\vec{X})=\sqrt{a^{2}+1}.
\end{equation*}
Combining this with (\ref{rnd}) gives that $\int_{0}^{a}g_{1}(\xi
)d\xi =ac_{a}$, which contradicts (\ref{oonk}) and so establishes
that $\theta _{a}<\arccos \frac{1}{w(0)}$.

From the past two paragraphs and the fact that $w$ is strictly
decreasing
from $w(0)$ to $1$ on $[0,1]$ we deduce that there exists a unique number $%
\xi _{a}\in (0,1)$ such that $\theta _{a}=\arccos \frac{1}{w(\xi
_{a})}$ and that $U=(\xi _{a},a)$. This in turn implies that
\begin{eqnarray*}
a &=&\int_{0}^{a}g_{0}(\xi )d\xi =\int_{0}^{a}\rho (\xi )\cos
\theta (\xi
)d\xi \\
&=&\int_{0}^{\xi _{a}}\frac{c_{a}}{\cos \theta (\xi )}\cos \theta
(\xi )d\xi
+\int_{\xi _{a}}^{a}c_{a}w(\xi )\cos \theta _{a}d\xi \\
&=&c_{a}\xi _{a}+c_{a}\cos \theta _{a}\left( E(\xi _{a},\alpha ;\vec{X}%
)-E(a,\alpha ;\vec{X})\right) \\
&=&c_{a}\xi _{a}+\frac{c_{a}}{w(\xi _{a})}\sqrt{(a-\xi
_{a})^{2}+(1-\xi
_{a})^{2}} \\
&=&c_{a}\xi _{a}+\frac{c_{a}}{a+1-2\xi _{a}}\left( (a-\xi
_{a})^{2}+(1-\xi
_{a})^{2}\right) \\
&=&\frac{c_{a}}{a+1-2\xi _{a}}\left( a\xi _{a}+\xi _{a}-2\xi
_{a}^{2}+a^{2}-2a\xi _{a}+\xi _{a}^{2}+1-2\xi _{a}+\xi _{a}{}^{2}\right) \\
&=&\frac{c_{a}}{a+1-2\xi _{a}}\left( a^{2}-a\xi _{a}+1-\xi
_{a}\right) .
\end{eqnarray*}

\smallskip So we have
\begin{equation}
c_{a}=\frac{a^{2}+a-2a\xi _{a}}{a^{2}-a\xi _{a}+1-\xi _{a}}
\label{first}
\end{equation}

We also have

\begin{eqnarray*}
1 &=&\int_{0}^{a}g_{1}(\xi )d\xi =\int_{0}^{a}\rho (\xi )\sin
\theta (\xi
)d\xi \\
&=&\int_{0}^{\xi _{a}}\frac{c_{a}}{\cos \theta (\xi )}\sin \theta
(\xi )d\xi
+\int_{\xi _{a}}^{a}c_{a}w(\xi )\sin \theta _{a}d\xi \\
&=&c_{a}\int_{0}^{\xi _{a}}\tan \theta (\xi )d\xi +c_{a}\sin
\theta
_{a}\left( E(\xi _{a},\alpha ;\vec{X})-E(a,\alpha ;\vec{X})\right) \\
&=&c_{a}\int_{0}^{\xi _{a}}\sqrt{w^{2}(\xi )-1}d\xi +c_{a}\sqrt{1-\frac{1}{%
w^{2}(\xi _{a})}}\sqrt{(a-\xi _{a})^{2}+(1-\xi _{a})^{2}}
\end{eqnarray*}
We have already calculated another expression for $w^{2}(\xi )-1$
in (\ref {gac}) and (\ref{agac}), so we can substitute it in both
terms of the preceding line to get
\begin{displaymath}
\begin{split}
1 &=c_{a}\int_{0}^{\xi _{a}}\sqrt{\frac{2(a-\xi )(1-\xi )}{(a-\xi
)^{2}+(1-\xi )^{2}}}d\xi+\\
&+c_{a}\frac{1}{w(\xi _{a})}\sqrt{\frac{2(a-\xi _{a})(1-\xi
_{a})}{(a-\xi _{a})^{2}+(1-\xi _{a})^{2}}}\sqrt{(a-\xi
_{a})^{2}+(1-\xi
_{a})^{2}}=\\
&=c_{a}\int_{0}^{\xi _{a}}\sqrt{\frac{2(a-\xi )(1-\xi )}{(a-\xi
)^{2}+(1-\xi )^{2}}}d\xi +c_{a}\sqrt{\frac{2(a-\xi _{a})(1-\xi
_{a})\left( (a-\xi _{a})^{2}+(1-\xi _{a})^{2}\right)
}{a+1-2\xi _{a}}}.\\
\end{split}
\end{displaymath}
\smallskip This latter formula can be rewritten as
\begin{equation}\label{vsmb}
\frac{1}{c_{a}}=\int_{0}^{\xi _{a}}\sqrt{\frac{2(a-\xi )(1-\xi
)}{(a-\xi )^{2}+(1-\xi )^{2}}}d\xi +\sqrt{\frac{2(a-\xi
_{a})(1-\xi _{a})}{a+1-2\xi _{a}}\cdot (a-\xi _{a})^{2}+(1-\xi
_{a})^{2}}.
\end{equation}

If we now substitute for $c_{a}$ in this equation, using (\ref
{first}) we will obtain a rather complicated equation for $\xi
_{a}$, which we will investigate further in the next section.

On a more simple level, we can use (\ref{first}) to obtain estimates for $%
c_{a}$ from above and below.
\begin{equation*}
\inf_{t\in (0,1)}\frac{a^{2}+a-2at}{a^{2}+1-(a+1)t}\le c_{a}\le
\sup_{t\in (0,1)}\frac{a^{2}+a-2at}{a^{2}+1-(a+1)t}
\end{equation*}
The function $t\mapsto \frac{a^{2}+a-2at}{a^{2}+1-(a+1)t}$ like
any function of the form $A\frac{b-t}{c-t}$ where $A$ , $b$ and
$c$ are positive constants, is either an increasing or decreasing
function on any interval which does not contain the point where
its denominator vanishes. In this case, its minimum on $[0,1]$
equals $1$ and is attained at $t=1$. Its maximum is
$\frac{a^{2}+a}{a^{2}+1}$ and is attained at $t=0$. The maximum
value of $\frac{a^{2}+a}{a^{2}+1}$ as $a$ ranges over $[1,\infty
)$ is
attained at $a=1+\sqrt{2}$ and is thus equal to $\frac{4+3\sqrt{2}}{4+2\sqrt{%
2}}$ which is approximately equal to $1.2071$.

\smallskip

\subsection{Some numerical experiments.}

\smallskip In this section we present some numerical experiments, which lead us to
a guess for the approximate value of the $K$-divisibility constant
of $(\ell_2^2,\ell^\infty_2)$, namely
$\gamma(\ell_2^2,\ell^\infty_2)\approx 1.0304$. Fix some value of
$a$ and try to find the corresponding value of $x=\xi _{a}$ by
defining
\begin{displaymath}
\begin{split}
f(x)&=\int_{0}^{x}\sqrt{\frac{2(a-t)(1-t)}{(a-t)^{2}+(1-t)^{2}}}%
dt+\\
&+\sqrt{\frac{2(a-x)(1-x)\left( (a-x)^{2}+(1-x)^{2}\right) }{a+1-2x}}-\frac{%
a^{2}-ax+1-x}{a^{2}+a-2ax}\\
\end{split}
\end{displaymath} and solving the equation (\ref{vsmb}) which is simply
$f(x)=0$. We are using ``Maple'' via its interface with ``Scientific
Workplace''. We will fix some values of $a$ and then try to find
$x\in (0,1)$ such that $f(x)=0$. We are currently ignoring the
question of whether such an $x$ is unique. To find the corresponding
value of $c_{a}$ we compute $g(x)=\frac{a^{2}+a-2ax}{a^{2}-ax+1-x}$.

Here is a table which summarizes some of our numerical
experiments, and which indicates that maybe the value of $\gamma$
is approximately $1.0304$:

\smallskip

\begin{center}
\begin{tabular}{ccc}
$a$ & $x$ & $c_{a}$ \\
&  &  \\
1.2 & $.94667221295$ & $1.\,0298$ \\
1.25 & $.94778089315$ & $1.\,0304$ \\
1.3 & $.94840470115$ & $1.\,0304$ \\
1.275 & $.94811047015$ & $1.\,0304$ \\
1.5 & $.95139101435$ & $1.\,0279$ \\
1.6 & $.95340037845$ & $1.\,0259$ \\
1.8 & $.95781371025$ & $1.\,0217$ \\
1.2 & $.94667221295$ & $1.\,0298$ \\
2 & $.96218058915$ & $1.\,0179$ \\
2.2 & $.96618489325$ & $1.\,0148$ \\
$1+\sqrt{2}$ & $.96997017725$ & $1.\,0121$ \\
3 & $.977870722252$ & $1.\,0073$ \\
&  &
\end{tabular}
\end{center}

\section{\protect\smallskip \label{appendix}Appendix: The couple $(\ell
_{n}^{2},\ell _{n}^{\infty })$ is an exact Calder\'{o}n couple
when $n=2$, but not when $n=8$.}

Suppose that $\vec{X}=(\ell _{2}^{2},\ell _{2}^{\infty })$ and that $%
f=(f_{0},f_{1})$ and $g=(g_{0},g_{1})$ are two points in
$\Bbb{R}^{2}$ which satisfy $K(t,g;\vec{X})\le K(t,f;\vec{X})$ for
all $t>0$. We will show that there exists an operator
$T:\vec{X}\rightarrow \vec{X}$ with norm $\left\| T\right\|
_{\vec{X}\rightarrow \vec{X}}\le 1$ such that $Tf=g$. We can of
course assume without loss of generality that $f_{0}\ge f_{1}\ge 0$ and $%
g_{0}\ge g_{1}\ge 0$. The $K$-functional inequalty satisfied by
$f$ and $g$ is equivalent to an $E$-functional inequality which
can be written as
$$\left( f_{0}-\min (t,f_{0})\right) ^{2}+\left(
f_{1}-\min (t,f_{1})\right) ^{2}\ge \left( g_{0}-\min
(t,g_{0})\right) ^{2}+\left( g_{1}-\min (t,g_{1})\right) ^{2}$$
and which holds for all $t>0$.

It is clear that $f_{0}\ge g_{0}$. (Otherwise we get a
contradiction by
choosing $t=(f_{0}+g_{0})/2$.) By setting $t=0$ we also have that $%
f_{0}^{2}+f_{1}^{2}\ge g_{0}^{2}+g_{1}^{2}$. This means that the
condition
\begin{equation}
\int_{0}^{t}f^{*}(s)^{2}ds\ge \int_{0}^{t}g^{*}(s)^{2}ds
\label{ih}
\end{equation}
holds for $t=0$, $1$ and for all $t\ge 2$. Since both sides of
(\ref{ih}) are affine functions on $[0,1]$ and $[1,2]$ it follows
that (\ref{ih}) holds for all $t>0$. Then we can apply the theorem
and proof of Lorentz and Shimogaki \cite{lorshi} to construct the
required operator $T$.

In contrast to the preceding calculation let us now show that
$(\ell
_{n}^{2},\ell _{n}^{\infty })$ is not an exact Calder\'{o}n couple for all $%
n\ge 8$. It is conceivable that a similar result also holds for
other smaller values of $n$. We recall that it was shown by Sparr
\cite{sparr} Example 5.1 that the five-dimensional version of the
dual couple $(\ell_5^{1},\ell_5^{2})$ is not an exact Calder\'{o}n
couple. (Note also that Brudnyi and Shteinberg \cite{bs} have
studied relations between the Calder\'{o}n constants for a
finite-dimensional couple and for its dual.)

We set $n=8$ and consider the two vectors $f$ and $g$ in
$\Bbb{R}^{8}$ given by $$f=(3,1,1,1,1,1,1,1)\text{ and }
g=(2,2,2,2,0,0,0,0).$$ Then $E(t,f;\ell ^{2},\ell ^{\infty })=\inf
\left\{ \left\| f-h\right\| _{\ell ^{2}}:h\in \ell ^{\infty
},\left\| h\right\| _{\ell ^{\infty }}\le t\right\} $ satisfies

\begin{equation*}
E(t,f;\ell ^{2},\ell ^{\infty })=\left\{
\begin{array}{lll}
\sqrt{(3-t)^{2}+7(1-t)^{2}} & , & 0\le t\le 1 \\
3-t & , & 1<t<3 \\
0 & , & t\ge 3
\end{array}
\right.
\end{equation*}
a\smallskip nd the corresponding error functional for $g$ is given
by

\begin{equation*}
E(t,g;\ell ^{2},\ell ^{\infty })=\left\{
\begin{array}{lll}
4-2t & , & 0\le t\le 2 \\
0 & , & t>2
\end{array}
\right. .
\end{equation*}
Clearly $E(t,g;\ell ^{2},\ell ^{\infty })\le E(t,f;\ell ^{2},\ell
^{\infty }) $ for all $t\ge 0$. So, if $(\ell _{8}^{2},\ell
_{8}^{\infty })$ is an exact Calder\'{o}n couple, there should be
a linear operator $T:(\ell _{8}^{2},\ell _{8}^{\infty
})\rightarrow (\ell _{8}^{2},\ell _{8}^{\infty })$ with $\left\|
T\right\| _{(\ell _{8}^{2},\ell _{8}^{\infty })\rightarrow (\ell
_{8}^{2},\ell _{8}^{\infty })}\le 1$ such that $Tf=g$. Suppose
that such a $T$ exists, and let $\lambda :\Bbb{R}^{8}\rightarrow
\Bbb{R}$ be the linear functional obtained by defining $\lambda
(h)=\frac 1 8 \sum_{i=1}^4(Th)_i$. Then $\lambda $ is given by the formula $%
\lambda (h)=\sum_{j=1}^{8}\lambda _{j}h_{j}$ where the numbers
$\lambda _{j}$ must satisfy
\begin{equation}
\sum_{j=1}^{8}|\lambda _{j}|\le 1  \label{lie}
\end{equation}
and also
\begin{equation}
2\sum_{j=1}^{8}\lambda _{j}h_{j}\le
\sqrt{\sum_{j=1}^{8}h_{j}^{2}}. \label{rpfgs}
\end{equation}
The condition $Tf=g$, i.e.\ $\lambda (f)=2$, implies that equality
holds in (\ref{rpfgs}) when $h=f$. By standard
facts about the Cauchy-Schwartz inequality, this in turn implies that $%
(\lambda _{1},\lambda _{2},...,\lambda _{8})=\frac{1}{8}f$. But
this contradicts (\ref{lie}) and so we have shown that $(\ell
_{8}^{2},\ell _{8}^{\infty })$ is not an exact Calder\'{o}n
couple.


\begin{thebibliography}{99}

\bibitem{a2}
Y.\ Ameur, A new proof of Donoghue's interpolation theorem,
Journal of Function Spaces and Applications\ 3 (2004), 253--265.

\bibitem{ameur}  Y.\ Ameur, The Calder\'{o}n problem for Hilbert couples,
Ark.\ Mat.\ 41 (2003), 203--231.

\bibitem{bs1}  C.\ Bennett and R.\ Sharpley, K-divisibility and a theorem of
Lorentz and Shimogaki, Proc.\ Amer.\ Math.\ Soc., 96 (1986), 585--592.

\bibitem{bsh}  C.\ Bennett and R.\ Sharpley, \textit{Interpolation of
Operators,} Academic Press, Boston 1988.

\bibitem{bl}  J.\ Bergh and J.\ L\"{o}fstr\"{o}m, \textit{Interpolation
spaces. An Introduction,} Springer, Berlin 1976.

\bibitem{bk1}  Ju.\ A.\ Brudny\v {i}, N.\ Ja.\ Krugljak, Real interpolation
functors, Dokl.\ Akad.\ Nauk SSSR, 256 (1981), 14--17 = Soviet Math.Dok. 23
(1981), 6--8.

\bibitem{bk}  Y.\ Brudnyi and N.\ Krugljak, \textit{Interpolation functors
and interpolation spaces, Volume 1,} North Holland, Amsterdam 1991.

\bibitem{bs} Y.\ Brudnyi and A. Shteinberg, \textit{Calder\'{o}n
constants of finite-dimensional couples,} Israel J. Math. 101
(1997), 289--322.

\bibitem{ckdiv}  M.\ Cwikel, $K$-divisibility of the $K$-functional and
Calder\'{o}n couples, Ark.Mat. 22 (1984), 39--62.

\bibitem{cw1}  M. Cwikel, The $K$-divisibility constant for couples of
Banach lattices. J.\ Approx.\ Th. 124 (2003) 124--136.

\bibitem{cjm}  M.\ Cwikel, B.\ Jawerth and M.\ Milman, On the fundamental
lemma of interpolation theory. J.\ Approx.\ Th.\ 60 (1990) 70--82.

\bibitem{ckeich}  M.\ Cwikel and U.\ Keich, Optimal decompositions for the $K
$-functional for a couple of Banach lattices. Arkiv f\"{o}r Matematik. 39
(2001) 27--64.

\bibitem{cwikoz}  M.\ Cwikel and I.\ Kozlov, Interpolation of weighted $L^{1}
$ spaces - a new proof of the Sedaev-Semenov theorem. Illinois J.\ Math. 46
(2002) 405--419.

\bibitem{cn}  M.\ Cwikel and P.\ G.\ Nilsson, Interpolation of weighted
Banach lattices, Memoirs Amer.\ Math.\ Soc. 165 (2003) no.\ 787,
1-105.

\bibitem{cp}  M.\ Cwikel, M. and J.\ Peetre, Abstract K and J spaces. J.\
Math.\ Pures et Appl.\ {60} (1981), 1--50.

\bibitem{daa}  A.\ A.\ Dmitriev, On the exactness of Peetre $K$%
-interpolation method. In: \textit{Operator methods in Differential Equations%
}, Voronezh, 1979, 32--40. (Russian).

\bibitem{do}
W.\ Donoghue, The interpolation of quadratic norms.\ Acta\ Math.\
{118} (1967), 251--270.

\bibitem{kr}  N.~Ja.\ Krugljak, On the K-divisibility constant of the couple
$(C,C^{1})$, in: \textit{Analysis of the Theory of Functions of Several Real
Variables}, Yaroslavl, 1981, 37--44. (Russian)

\bibitem{lorshi}  G.\ G.\ Lorentz and T.\ Shimogaki, Interpolation theorems
for the pairs of spaces $(L^{p},L^{\infty })$ and $(L^{1},L^{q})$. Trans.\
Amer.\ Math.\ Soc., 159 (1971), 207--222.

\bibitem{po}  T.\ S.\ Podogova, On a property of modulus of continuity, in:
\textit{Analysis of the Theory of Functions of Several Real Variables},
Yaroslavl, 1982, 84--89. (Russian)

\bibitem{shv}  P.\ Shvartsman, A geometrical approach to the K-divisibility
problem, Israel J.\ Math., 103 (1998), 289--318.

\bibitem{sed}  A.\ A.\ Sedaev, Description of interpolation spaces for the
couple $(L_{\alpha _{0}}^{p},L_{\alpha _{1}}^{p})$ and some
related problems. { Dokl.\ Akad.\ Nauk SSSR,} { 209} (1973),
799-800 (Russian), { Soviet Math.\ Dokl.,} { 14} (1973), 538-541.

\bibitem{sedsem}  A.\ A.\ Sedaev and E.\ M.\ Semenov, On the possibility of
describing interpolation spaces in terms of Peetre's $K$-method. \textit{%
Optimizaciya,} \textbf{4} (1971), 98--114 (Russian).

\bibitem{sparr}  G.\ Sparr, Interpolation of weighted $L^{p}$ spaces.
\textit{Studia Math.,} \textbf{62} (1978), 229--271.

\bibitem{sw}  E.\ M.\ Stein and G.\ Weiss, Interpolation of operators with
change of measures, Trans.\ Amer.\ Math.\ Soc.\ 87 (1958), 159--172.
\end{thebibliography}
\end{document}